\documentclass[12pt]{amsart}
\usepackage{amsmath}
\usepackage{amssymb}
\usepackage{graphicx}
\usepackage{color}
\usepackage{tcolorbox}
\usepackage{comment}

\definecolor{violet}{rgb}{0.5,0,0.5}
\definecolor{orange}{cmyk}{0,0.3,0.7,0}

\newcommand{\RR}{\mathbb{R}}

\newcommand{\eqdef}{\stackrel{{\mathrm {def}}}{=}}

\renewcommand{\colon}{:\,}

\numberwithin{equation}{section}
\newtheorem{theorem}{Theorem}[section]
\newtheorem{proposition}[theorem]{Proposition}
\newtheorem{lemma}[theorem]{Lemma}
\newtheorem{corollary}[theorem]{Corollary}
\newtheorem{remark}[theorem]{Remark}
\newtheorem{example}[theorem]{Example}

\allowdisplaybreaks

\newcommand{\Square}{$\sqcap$\hskip -1.5ex $\sqcup$}

\begin{document}

\title[Strong Comparison principle for Degenerate Elliptic Problems]%
{ On the Strong Comparison Principle for\\
  Degenerate Elliptic Problems with Convection}

\author[J.~Benedikt]{Ji\v{r}\'i Benedikt}
\author[P.~Girg]{Petr Girg}
\author[L.~Kotrla]{Luk\'{a}\v{s} Kotrla}
\author[P.~Tak\'a\v{c}]{Peter Tak\'a\v{c}}

\address{Ji\v{r}\'{\i} Benedikt, Petr Girg, Luk\'a\v{s} Kotrla\newline
         Department of Mathematics and NTIS,
         Faculty of Applied Scences,\hfil\break
         University of West Bohemia, Univerzitn\'{\i} 8,
         CZ--301\,00~Plze\v{n}, Czech Republic}
\email{benedikt@kma.zcu.cz, pgirg@kma.zcu.cz, kotrla@kma.zcu.cz}

\address{Peter Tak\'a\v{c}\newline
         Institut f{\"u}r Mathematik, Universit{\"a}t Rostock,
         Ulmenstra{\ss}e 69, D--18055 Rostock, Germany}
\email{peter.takac@uni-rostock.de}



\begin{abstract}
The weak and strong comparison principles
({\em\bfseries WCP\/} and {\em\bfseries SCP\/}, respectively)
are investigated for quasilinear elliptic boundary value problems
with the $p$\--Laplacian in one space dimension,
\begin{math}
  \Delta_p(u)\eqdef
  \frac{\mathrm{d}}{\mathrm{d}x} \left( |u'|^{p-2} u'\right) \,.
\end{math}
We treat the ``degenerate'' case of $2 < p < \infty$
and allow also for the nontrivial
{\it\bfseries convection velocity\/}
$b\colon [-1,1]\to \RR$ in the underlying domain $\Omega = (-1,1)$.
We establish the {\em\bfseries WCP\/}
under a rather general, ``natural sufficient condition''
on the convection velocity, $b(x)$, and
the reaction function, $\varphi(x,u)$.
Furthermore, we establish also the {\em\bfseries SCP\/}
under a number of various additional hypotheses.
In contrast, with these hypotheses being violated,
we construct also a few rather natural counterexamples
to the {\em\bfseries SCP\/} and discuss their applications to
an interesting classical problem of fluid flow in porous medium,
{\em ``seepage flow of fluids in inclined bed''\/}.
Our methods are based on a mixture of classical and new techniques.
\end{abstract}

\maketitle

\section{Introduction}
\label{s:intro}

The {\em\bfseries weak\/} and
{\em\bfseries strong comparison principles\/}
({\em\bfseries WCP\/} and {\em\bfseries SCP\/}, respectively,
 for brevity)
play an important role in the theory of elliptic and parabolic problems
for partial differential equations.
In our present work we focus on such problems in one space dimension
that ``idealize'' long, thin domains
$\Omega\subset \mathbb{R}^N$ of type
$\Omega = \omega\times (a,b)$
where solutions are expected to be rather homogeneous
on the relatively small cross\--section
$\omega\subset \mathbb{R}^{N-1}$ (a domain in $\mathbb{R}^{N-1}$)
of the relatively long domain $\Omega = \omega\times (a,b)$.
Here, we assume that
$\mathrm{diameter}(\omega) \ll |b-a|$ meaning that the ratio
$\mathrm{diameter}(\omega) / |b-a|$ is a sufficiently small number
${}< 1$.

In general, in an arbitrary domain
$\Omega\subset \mathbb{R}^N$ of space dimension $N$ ($N\geq 1$),
we consider the elliptic partial differential equation
\begin{equation}
\label{e:ell_PDE}
\left\{\hspace{2.00mm}
\begin{aligned}
{}- \Delta_p u  - \mathbf{b}(x)\cdot \nabla u + \varphi(x,u)
  & = f(x) &&\qquad \mbox{ for }\, x\in \Omega \,;
\\
  u & = 0  &&\qquad \mbox{ on }\, \partial\Omega \,.
\end{aligned} 
\right.\hfill
\end{equation}
Here, we assume that $2 < p < \infty$,
\begin{math}
  \Delta_p(u)\eqdef \mathrm{div}\left( |\nabla u|^{p-2} \nabla u\right)
\end{math}
stands for $p$\--Laplacian,
$\mathbf{b}\colon \Omega\to \mathbb{R}^N$
is an essentially bounded measurable vector field
which stands for the {\it\bfseries convection velocity\/},
often termed also
{\it\bfseries drift or advection velocity\/},
the {\it\bfseries reaction function\/}
$\varphi\colon \Omega\times \mathbb{R}\to \mathbb{R}$
is continuously partially differentiable with respect
to the second variable, $u$, with the partial derivative
$\frac{\partial\varphi}{\partial u} (x,u)$ $> 0$ for all
$(x,u)\in \Omega\times (\mathbb{R}\setminus \{ 0\})$, and
$\varphi(x,0) = 0$ for all $x\in \Omega$.
Finally, $f\colon \Omega\to \mathbb{R}$ describes the volumetric source
density, $f\in L^{\infty}(\Omega)$; sometimes we will take
$f\in L^{\infty}(\Omega)\cap C(\Omega)$ or similar.

Under some reasonable hypotheses on the vector field $\mathbf{b}$
and the function $\varphi$,
a standard \emph{weak comparison} result for
this Dirichlet boundary value problem is established in the monograph by
{\sc D.\ Gilbarg} and {\sc N.~S.\ Trudinger}
\cite[Sect.~10, pp.\ 263--271]{GilbargTrudinger2001},
Theorem 10.1 (p.~263) and Theorem 10.7 (p.~268).
In contrast, the corresponding \emph{strong comparison} and
\emph{Hopf\--type comparison} results are, in general,
{\bf not valid\/} under analogous hypotheses; cf.\
{\sc M.\ Guedda} and {\sc L.\ V\'eron} \cite{GueddaVeron},
Prop.\ 2.1 (pp.\ 886--887) and Prop.\ 2.2 (p.~888),
and later work in
{\sc M.\ Cuesta} and {\sc P.\ Tak\'{a}\v{c}}
\cite{CuestaTakac1998, CuestaTakac2000}. 
We begin this article by giving a few simple counterexamples
to the strong comparison principle. 
Our counterexamples will be constructed in the one\--dimensional domain
$\Omega = (-1,1)\subset \mathbb{R}$.
In fact, they have motivated the entire work presented in this article.
Our main goal is to find reasonable sufficient conditions
on the (now) scalar field
$\mathbf{b}\equiv b\colon [-1,1]\to \mathbb{R}^1$
and the reaction function $\varphi$
that would exclude such counterexamples;
we will formulate sufficient conditions on $b$ and $\varphi$
that will guarantee the validity of the strong comparison principle.
While assuming $\varphi\equiv 0$, we formulate sufficient conditions on
the {\it convection velocity\/}
$b$ in Section~\ref{s:strong_critical}.
The main purpose of the counterexamples will be to ``illustrate''
the necessity of our sufficient conditions.

\subsection{Mathematical Problems of weak and strong
            Elliptic Comparison Principles}
\label{ss:intro-math}

The one\--dimensional analogue of
the elliptic partial differential equation \eqref{e:ell_PDE}
has been investigated in the work by
{\sc M.\ Cuesta} and {\sc P.\ Tak\'{a}\v{c}}
\cite{CuestaTakac1998, CuestaTakac2000}.
More precisely, this work treats the special case with
$b(x)\equiv 0$ in an open interval
$\Omega = (-1,1)\subset \mathbb{R}^1$ and
$\varphi(x,u) = \lambda\, |u|^{p-2} u$ for all
$x\in (-1,1)$ and all $u\in \RR = (-\infty,\infty)$,
with a constant $\lambda\in (0,\infty)$ large enough,
say, $\lambda\geq \lambda_p > 0$.
In this way, the original problem~\eqref{e:ell_PDE} is reduced to
the following two\--point boundary value problem,
\begin{equation}
\label{new:ell_PDE}
\left\{\hspace{2.00mm}
\begin{aligned}
{}- \left( |u'|^{p-2} u'\right)' + \lambda\, |u|^{p-2} u
& = f(x) \qquad \mbox{ for }\, x\in (-1,1) \,;
\\
  u(\pm 1) & = 0 \,.
\end{aligned} 
\right.
\end{equation}
This clearly is a variational problem on the Sobolev space
$W_0^{1,p}(-1,1)$ with the left\--hand side
\begin{equation*}
  [\mathcal{A}u](x)\eqdef
{}- \left( |u'|^{p-2} u'\right)' + \lambda\, |u|^{p-2} u \,,\quad
    x\in (-1,1) \,,
\end{equation*}
being a {\it monotone operator\/}
$\mathcal{A}\colon W_0^{1,p}(-1,1)\to W^{-1,p'}(-1,1)$
from its domain $W_0^{1,p}(-1,1)$ into the dual space
\begin{math}
  W^{-1,p'}(-1,1) = \left[ W_0^{1,p}(-1,1) \right]'
\end{math}
where $p'= p/(p-1)$.
We note the continuous imbedding
\begin{math}
  L^{\infty}(-1,1) = \left[ L^1(-1,1) \right]'
  \hookrightarrow W^{-1,p'}(-1,1)
\end{math}
that is dual to the Sobolev imbedding
\begin{math}
  W_0^{1,p}(-1,1)\hookrightarrow L^1(-1,1) \,.
\end{math}
The strict monotonicity of $\mathcal{A}$ hinges upon the fact that
\begin{math}
    \frac{\partial\varphi}{\partial u} (x,u)
  = \lambda\, (p-1)\, |u|^{p-2} > 0
\end{math}
is valid for all
$(x,u)\in \Omega\times (\mathbb{R}\setminus \{ 0\})$.
It forces the {\it\bfseries WCP\/}
({\it weak comparison principle\/},
 proved in Proposition~\ref{prop-Weak_Comp}),
by a standard variational argument that we will explain in
Section~\ref{s:weak_comp} below.
In this step, what matters is the following observation:
\begin{quote}
\begingroup\sl
The larger the constant $\lambda\geq 0$,
the ``stronger'' the strict monotonicity of $\mathcal{A}$.
\endgroup
\end{quote}
For $\lambda = 0$, even the {\it\bfseries SCP\/}
({\it strong comparison principle\/})
for problem~\eqref{new:ell_PDE} is established in
\cite{CuestaTakac1998, CuestaTakac2000}.
Nevertheless, the strong comparison principle fails to hold for every
$\lambda > 0$ large enough, say, $\lambda\geq \lambda_p > 0$,
as is shown in
{\sc M.\ Cuesta} and {\sc P.\ Tak\'{a}\v{c}}
\cite[Sect.~4, Example 4.1, pp.\ 740--741]{CuestaTakac2000};
see Example~\ref{exam-counter} in {\S}\ref{ss:weak-counter} below.

This is precisely the problem we treat in the present article:
\begingroup\sl
Conditions that are favorable for the validity of
the {\it\bfseries WCP\/}
may have a {\em negative effect\/} on the validity of
the {\it\bfseries SCP\/}.
\endgroup
The latter may fail to hold under those conditions
(in Example~\ref{exam-counter}).

Our strategy will be as follows:
In Section~\ref{s:weak_comp} we establish
the {\it\bfseries WCP\/}
under rather simple, but still ``reasonably'' general conditions
on the {\it\bfseries convection velocity\/}
$b\colon \Omega\to \mathbb{R}$ in the interval
$\Omega = (-1,1)\subset \mathbb{R}$
and the {\it\bfseries reaction function\/}
$\varphi\colon \Omega\times \mathbb{R}\to \mathbb{R}$.

In the remaining parts of this work we then {\bf assume} the validity of
the {\it\bfseries WCP\/}
and try to formulate ``reasonable'' sufficient conditions
that guarantee the validity of the {\it\bfseries SCP\/}.
We also investigate the ``sharpness'' of our sufficient conditions
in Section~\ref{s:strong_critical}.
We provide an interesting new counterexample to the {\it\bfseries SCP\/}
in Example~\ref{exam-SCP-Girg} for the special case $\varphi\equiv 0$.

A classical treatment of
the {\it\bfseries WCP\/} and {\it\bfseries SCP\/}
can be found in the monographs by
{\sc P.\ Pucci} and {\sc J.\ Serrin} \cite{Pucci-Serr_book}
(see also the article \cite{Pucci-Serrin_2007})
and
{\sc M.~H.\ Protter} and {\sc H.~F.\ Weinberger}
\cite{Protter-Weinberg}.
Both these books use the {\em\bfseries Hopf boundary point lemma\/},
due to {\sc E.~Hopf}, in order to derive
{\it\bfseries SCP\/} from {\it\bfseries WCP\/}, see e.g.\
\cite[Chapt~2, Sect.~3]{Protter-Weinberg}, Theorem~7 on p.~65.
We will use a modification of this lemma adapted to the $p$\--Laplacian
in the works by
{\sc P.\ Pucci} and {\sc J.\ Serrin}
\cite[{\S}5.4]{Pucci-Serr_book}, Theorem 5.5.1 on p.~120; cf.\ also
\cite[{\S}5.1]{Pucci-Serr_book}, Example on p.~104,
{\sc P.\ Tolksdorf}
\cite[Prop.\ 3.2.1 and 3.2.2, p.~801]{Tolksdorf-1983}, and
{\sc J.~L.\ V\'azquez} \cite[Theorem~5, p.~200]{Vazquez-1984}.
An important role in our application of {\sc Hopf}'s lemma is played by
the result ({\sc Hopf}--{\sc Oleinik} lemma)
established recently in the article by
{\sc D.~E.\ Apushkinskaya} and {\sc A.~I.\ Nazarov}
\cite[Theorem 2.1, p.~681]{Apush-Nazarov-2019}.
Last but not least, we should mention that the first sharp results
on the {\it\bfseries SCP\/} for the $p$\--Laplacian, $\Delta_p$,
were obtained in
{\sc M.\ Guedda} and {\sc L.\ V\'eron} \cite{GueddaVeron},
Prop.\ 2.1 (pp.\ 886--887) and Prop.\ 2.2 (p.~888).

\par\vskip 10pt
Several engineering problems connected with
the quasilinear ``diffusion'' operator $p$\--Laplacian, $\Delta_p$,
have been mentioned and sketched in the recent article
by the present authors,
{\sc J.\ Benedikt}, {\sc P.\ Girg}, {\sc L.\ Kotrla}, and
{\sc P.\ Tak\'a\v{c}} \cite{Benedikt-Kotrla2018}.
The two primary interests are focused on
{\sl\bfseries (i)}
{\sl water supplies for larger population communities\/}
and
{\sl\bfseries (ii)}
{\sl oil exploration and polution\/}.
For interesting applications of our results in the real world,
we refer the reader to the recent work by
{\sc R.~K.\ Bansal} \cite{Bansal2017}
for engineering\--type research on these topics.
There, a mathematical model of
{\em ``seepage flow of fluids in inclined bed''\/}
is presented and investigated.
Problems under consideration in our paper can be thought of being
stationary cases of time\--dependent problems studied
in~\cite{Bansal2017}.
In our case we consider the power law
(with $p>2$ corresponding to ``slow'' diffusion, see
 Paragraph \ref{ss:intro-fluid} below)
instead of Darcy's law; cf.\
\cite[Sect.~2]{Benedikt-Kotrla2018}
for ``historical'' justification and numerous references.
Indeed, it has been observed experimentally for a very slow flow
(in terms of Reynolds' number) in a porous medium,
that Darcy's law should be replaced by another law, namely,
the power law corresponding to the case $p>2$.
We refer the reader to the article by
{\sc J.~P.\ Soni}, {\sc N.\ Islam}, and {\sc P.\ Basak}
\cite[p.~239]{SoIsBa}
for a nice overview of the values of the constant $p$
for various materials where $p$ ranges from $2.12$ to $4.74$.
Note that the constant $n = p-1$ is used in \cite{SoIsBa} in place of~$p$.

\par\vskip 10pt
As far as applications of the two comparison principles,
{\it\bfseries WCP\/} and {\it\bfseries SCP\/},
to these engineering problems are concerned,
the following difference between them is of current interest:
\begin{quote}
If the source term $f(x)$ on the right\--hand side
in problem~\eqref{e:ell_PDE} is increased to a function $g(x)$
($g\geq f$ almost everywhere (a.e.) in $\Omega$)
only {\em\bfseries locally\/} near a point $x_0\in \Omega$, i.e.,
$g(x) > f(x)$ holds for almost all $x\in B_r(x_0)\subset \Omega$
in a small ball\hfil\break
\begin{math}
  B_r(x_0)\eqdef \{ x\in \mathbb{R}^N\colon |x-x_0| < r\}
          \subset \Omega
\end{math}
of radius $r > 0$ centered at $x_0\in \Omega$,
\hfil\break
{\bf does\/} the solution $u(x)$ {\bf increase\/} to $v(x)$
\hfil\break
(for notation, see problems \eqref{eq:WCP} in
 the next paragraph, {\S}\ref{ss:weak-counter})
\begin{itemize}
\item[{\bf (i)}]
only {\em\bfseries locally\/} near the point $x_0\in \Omega$, i.e.,
$v(x) > u(x)$ holds for all $x\in \Omega$ near the point $x_0$, say,
for all $x\in B_R(x_0)\subset \Omega$,
but $v(x) = u(x)$ holds for all
$x\in \Omega\setminus B_R(x_0)$, with some radius $R$, $r < R < \infty$,
{\bf or}
\item[{\bf (ii)}]
the inequality $v(x) > u(x)$ holds for all $x\in \Omega$, that is,
{\em\bfseries globally\/} throughout the entire domain $\Omega$ ?
\end{itemize}
The first statement, {\bf (i)}, means that the {\it\bfseries SCP\/}
fails to hold, whereas the second statement, {\bf (ii)}, is equivalent to
the {\it\bfseries SCP\/}.
\end{quote}
%

\subsection{A Simple Counterexample to the Strong Comparison Principle}
\label{ss:weak-counter}

We begin by recalling a special case from the work by
{\sc M.\ Cuesta} and {\sc P.\ Tak\'{a}\v{c}}
\cite[Sect.~4, Example 4.1, pp.\ 740--741]{CuestaTakac2000},
where a counterexample to the SCP for
the two\--point boundary value problem in eq.~\eqref{new:ell_PDE} above
is constructed by making the special choice $b(x)\equiv 0$ and
$\varphi(x,u) = \lambda\, |u|^{p-2} u$ for all
$x\in (-1,1)$ and all $u\in \RR = (-\infty,\infty)$
in the following two\--point boundary value problems
with $2 < p < \infty$:
\begin{equation}
\label{eq:WCP}
\left\{\quad
\begin{aligned}
{}- \left( |u'|^{p-2} u'\right)' - b(x)\, u' + \varphi(x,u)
& = f(x) \,,
\\
{}- \left( |v'|^{p-2} v'\right)' - b(x)\, v' + \varphi(x,v)
& = g(x) \,;
\\
  \quad\mbox{ with }\quad
  u(\pm 1)\leq v(\pm 1) \quad\mbox{ and }\quad
  f(x) &\leq g(x) \;\mbox{ for }\, x\in (-1,1) \,.
\end{aligned}
\right.
\end{equation}
The constant $\lambda\in (0,\infty)$ that appears in the reaction function
$\varphi(x,u) = \lambda\, u^{p-2} u$
in eq.~\eqref{eq:WCP}
has to be chosen large enough, say, $\lambda\geq \lambda_p > 0$.
As we work with nonnegative solutions
$u,v\colon [-1,1]\to \mathbb{R}$ to problems \eqref{eq:WCP} only,
it suffices to define the reaction function
$\varphi(x,u) = \lambda\, u^{p-1}$ for all
$x\in (-1,1)$ and all $u\in \RR_+ = [0,\infty)$ only.

As the boundary value problems in \eqref{eq:WCP}
in our present article allow for the (additional) convection terms
$b(x)\, u'$ and $b(x)\, v'$, respectively,
our counterexample below is a slight modification of
that one mentioned above \cite[Example 4.1]{CuestaTakac2000}
(with $b(x)\equiv 0$).

\par\vskip 10pt
\begin{example}[{\rm A Counterexample to the} {\it\bfseries SCP\/}]
\label{exam-counter}
\begingroup\rm
Let us consider the two\--point boundary value problems \eqref{eq:WCP}
that have arised from the following homogeneous Dirichlet
boundary value problem:
\begin{equation}
\label{eq:intro:counter}
\left\{\hspace{2.00mm}
\begin{aligned}
{}- \left( |u'|^{p-2} u'\right)' - b(x)\, u' + \varphi(u)
& = f(x) \qquad \mbox{ for }\, x\in (-1,1) \,;
\\
  u(\pm 1) & = 0 \,.
\end{aligned} 
\right.
\end{equation}
Let the real numbers $p, \theta\in \RR$ be given with
$p > 2$ and $\theta > 1$ arbitrary, and
let us define the following functions, while still allowing
$b\colon [-1,1]\to \mathbb{R}$ to be arbitrary, but
continuously differentiable:
\begin{align}
\label{e:u_theta}
& u_{\theta}(x) \stackrel{\rm def}{=}
  1 - |x|^{\theta} \qquad\mbox{ for }\, x\in [-1,1] \qquad\mbox{ and }
\\
&
\label{e:f_theta}
\begin{aligned}
  f_{\theta}(x) &\stackrel{\rm def}{=}
{}- \left( |u_{\theta}'|^{p-2} u_{\theta}'\right)'
  - b(x)\, u_{\theta}' + \varphi(u_{\theta})
\\
& = (p-1)\theta^{p-1} (\theta - 1)\, |x|^{(\theta - 1)(p-1) - 1} 
  + \theta\, b(x)\cdot \mathrm{sgn}(x)\, |x|^{\theta - 1} 
  + \varphi\left( 1 - |x|^{\theta}\right) \,.
\end{aligned}
\end{align}

We notice that for $1 < \theta_1 < \theta_2 < \infty$ we have
$u_{\theta_1}(x) < u_{\theta_2}(x)$ for all
$x\in (-1,1) \setminus \{ 0\}$, whereas
$u_{\theta_1}(0) = u_{\theta_2}(0) = 1$.
In formula~\eqref{e:b(x)} below we make a suitable choice of
the {\it\bfseries convection velocity\/}
$b(x)$, take $\varphi\colon \RR\to \RR$ linear,
$\varphi(u) = \lambda\, u$ in formula~\eqref{e:varphi(x)},
and state sufficient conditions on the exponents
$\theta_1 < \theta_2$ and the constant $\lambda > 0$,
such that the following inequality
\begin{equation}
\label{eq:intro:comparison}
  f_{\theta_1}(x) < f_{\theta_2}(x)
    \quad\mbox{ holds for all }\, x\in (-1,1)\setminus \{ 0\} \,.
\end{equation} 
Consequently, the strong comparison principle for problem
\eqref{eq:intro:counter} fails to hold at the point $x=0$ owing to
$u_{\theta_1}(0) = u_{\theta_2}(0)$.

When looking for a convection velocity $b(x)$ with some desired properties,
it is worth of noticing that
the {\em\bfseries reflection transformation\/}
\begin{math}
  \widetilde{\phantom{x}}\colon [-1,1]\to [-1,1]
                        \colon x\mapsto {\tilde x}\eqdef {}- x
\end{math}
brings problem~\eqref{eq:intro:counter}
into the following equivalent ``reflection'' form:
\begin{equation}
\nonumber
\tag{\ref{eq:intro:counter}$\widetilde{\phantom{x}}$}
\left\{\hspace{2.00mm}
\begin{aligned}
{}- \left( |{\tilde u}'|^{p-2} {\tilde u}'\right)'
  + b(- {\tilde x})\, {\tilde u}' + \varphi({\tilde u})
& = f(- {\tilde x}) \qquad \mbox{ for }\, {\tilde x}\in (-1,1) \,;
\\
  {}- {\tilde u}(\pm 1) & = 0 \,.
\end{aligned} 
\right.
\end{equation}
Here,
${\tilde u}({\tilde x}) = u(- {\tilde x}) = u(x)$
is the new unknown function of the variable
${\tilde x} = {}- x$, so that its derivative equals to
\begin{math}
    {\tilde u}'({\tilde x})
  = \frac{ \mathrm{d}{\tilde u} }{ \mathrm{d}{\tilde x} } ({\tilde x})
  = \frac{\mathrm{d}}{ \mathrm{d}{\tilde x} }\, u(- {\tilde x})
  = {}- \frac{\mathrm{d}u}{ \mathrm{d}x } (- {\tilde x})
  = {}- u'(- {\tilde x}) = {}- u'(x) \,,
\end{math}
by the chain rule.
In particular, the convection velocity $b(x)$ is transformed into
the new convection velocity 
\begin{math}
  {\tilde b}({\tilde x}) = {}- b(- {\tilde x}) = {}- b(x)
\end{math}
in the new problem
(\ref{eq:intro:counter}$\widetilde{\phantom{x}}$)
above, with the derivative
${\tilde b}^{\prime}({\tilde x}) = b'(x)$.
Consequently, the sign of the convection velocity $b(x)$,
in the expression $\pm b(x)$,
does not matter in our counterexample.
However, a sign change in $\pm b(x)$ may cause a change of the constant
$\lambda > 0$ in our choice of the reaction function
$\varphi(u) = \lambda\, u$; cf.\
eq.~\eqref{e:varphi(x)} just below combined with the calculations in
eq.~\eqref{eq:varphi(x)}.

We verify inequality \eqref{eq:intro:comparison} by showing that
for every $\theta\in (\theta_1, \theta_2)$ we have
\begin{equation}
\label{e:df/d_theta>0}
  \frac{\partial f_{\theta}}{\partial \theta}(x) > 0
    \quad\mbox{ for all }\, x\in (-1,1)\setminus \{ 0\} \,;
\end{equation} 
cf.\ inequalities \eqref{ineq:b(x)/x}.
A reader who is interested in detailed calculations is referred to
\cite[Example 4.1]{CuestaTakac2000}
and to our calculations below leading to inequalities
\eqref{ineq:b(x)/x} and \eqref{cond:b(x)/x}.
\endgroup
\end{example}
\par\vskip 10pt

As a simple {\bf example}, owing to ineq.~\eqref{ineq:b(x)/x} below,
we may take

\par\vskip 10pt
\par\noindent
{\bf Example~\ref{exam-counter}} \underline{\rm continued:}$\quad$
\begingroup\rm
Let us take
\begin{align}
\label{e:p,theta}
& 2 < p < \infty \,,\quad
  \frac{p}{p-2} < \theta_1 < \theta_2 < \infty \,,
\\
\label{e:b(x)}
&
\begin{aligned}
  b(x)
&     \eqdef (p-1)^2 (\theta_1 - 1) \theta_1^{p-2}\,
             |x|^{ (p-2) (\theta_2 - 1) - 2 }\, x
\\
&     \equiv (p-1)^2 (\theta_1 - 1) \theta_1^{p-2}\,
             |x|^{ (p-2) (\theta_2 - 1) - 1 }\cdot \mathrm{sgn}(x)
      \quad\mbox{ for }\, x\in [-1,1] \,,
\end{aligned}
\\
\label{e:varphi(x)}
& \mbox{ and }\qquad
\left\{\hspace{2.00mm}
\begin{aligned}
& \varphi(x,s)\equiv \varphi(s)\eqdef \lambda\, s \quad\mbox{ for }\,
         (x,s)\in [-1,1]\times \mathbb{R}
\\
& \quad\mbox{ with the constant }\quad
  \lambda\geq 2 (p-1)^2\, (\theta_2 - 1) \theta_2^{p-1} > 0 \,.
\end{aligned} 
\right.
\end{align}
We observe that for every $x\in [-1,1]$ we have
\begin{equation}
\label{e:b'(x)}
  b'(x) = (p-1)^2 (\theta_1 - 1) \theta_1^{p-2}\,
          [ (p-2) (\theta_2 - 1) - 1 ]\,
          |x|^{ (p-2) (\theta_2 - 1) - 2 }\geq 0 \,.
\end{equation}
\endgroup
%
\par\vskip 10pt

We note that $\mathrm{sgn}(0) = 0$,
$\mathrm{sgn}(x) = 1$ if $x > 0$, and
$\mathrm{sgn}(x) = {}- 1$ if $x < 0$.

\par\vskip 10pt
{\it\bfseries Proof of\/} Ineq.~\eqref{e:df/d_theta>0}.
For $u(x)\equiv u_{\theta}(x)$ with $1 < \theta < \infty$
we make use of eq.~\eqref{e:u_theta} to calculate
\begin{align}
\nonumber
& u'(x)\equiv u_{\theta}'(x)
  = \frac{\mathrm{d}}{\mathrm{d}x} \left( 1 - |x|^{\theta}\right)
  = {}- \theta\, |x|^{\theta - 2}\, x
  = {}- \theta\, |x|^{\theta - 1}\cdot \mathrm{sgn}(x) \,,
\\
\label{eq:du/d_theta}
& \frac{\partial u_{\theta}}{\partial\theta}(x)
  = {}- |x|^{\theta}\, \log|x|\geq 0
  \quad\mbox{ for all }\, x\in [-1,1]\setminus \{ 0\} \,,
    \quad\mbox{ and }\quad
  \frac{\partial u_{\theta}}{\partial\theta}(0) = 0 \,.
\end{align}
It follows that
\begin{align*}
  |u'(x)|^{p-2} u'(x) & =
{}- \theta^{p-1}\, |x|^{(\theta - 1)(p-2) + \theta - 2} x =
{}- \theta^{p-1}\, |x|^{(\theta - 1)(p-1) - 1} x \,,
\\
{} - \left( |u'|^{p-2} u'\right)'(x) & \equiv
{} - \genfrac{}{}{}1{\mathrm{d}}{\mathrm{d}x}
     \left( |u'|^{p-2} u'\right)(x)
   = \theta^{p-1} (\theta - 1) (p-1) |x|^{(\theta - 1)(p-1) - 1} \,.
\end{align*}
Similarly, we use eq.~\eqref{e:f_theta} to get
\begin{align}
\label{eq:df/d_theta}
\begin{aligned}
  \frac{\partial f_{\theta}}{\partial\theta}(x)
& {}= (p-1)
    \left[ p\theta^{p-1} - (p-1)\theta^{p-2}\right]
    |x|^{(\theta - 1)(p-1) - 1}
\\
& {}+ (p-1)^2 \theta^{p-1}(\theta - 1)
      |x|^{(\theta - 1)(p-1) - 1}\, \log|x|
\\
& {}+ b(x)\cdot \mathrm{sgn}(x)\, |x|^{\theta - 1}
    + \theta\, b(x)\cdot{\rm sgn}(x)\, |x|^{\theta - 1}\, \log|x|
\\
& {}- \varphi'\left( 1 - |x|^{\theta}\right)\, |x|^{\theta}\, \log|x| \,.
\end{aligned}
\end{align}
In order to determine the sign of
${\partial f_{\theta}} / {\partial\theta}$,
we take into account that
$\varphi(x,s)\equiv \varphi(s) = \lambda\, s$,
by eq.~\eqref{e:varphi(x)}, and thus simplify
the right\--hand side of eq.~\eqref{eq:df/d_theta} as follows:
\begin{align}
\label{e:df/d_theta}
  \frac{\partial f_{\theta}}{\partial\theta}(x)
& {}= |x|^{\theta - 1}
    \left\{
    (p-1)\left[ p\theta^{p-1} - (p-1)\theta^{p-2}\right]\,
    |x|^{(\theta - 1)(p-2) - 1} + b(x)\cdot \mathrm{sgn}(x)
    \right\}
\\
\nonumber
&
\begin{aligned}
  {}+ |x|^{\theta - 1}\, \log|x|\cdot
&   \left\{ (p-1)^2\theta^{p-1}(\theta - 1)\, |x|^{(\theta - 1)(p-2) - 1}
\right.
  {}+ \theta\, b(x)\cdot \mathrm{sgn}(x)
\\
&
\left.
  {}- \lambda\, |x|
    \right\}
    \qquad\mbox{ for all }\, x\in [-1,1]\setminus \{ 0\} \,.
\end{aligned}
\end{align}
We observe that the desired inequality \eqref{eq:intro:comparison}
is satisfied provided the two expressions in braces
on the right\--hand side of equation \eqref{e:df/d_theta} above
have the correct sign, whenever $0 < |x| < 1$, that is to say,
\begin{align*}
& (p-1)\left[ p\theta^{p-1} - (p-1)\theta^{p-2}\right]\,
    |x|^{(\theta - 1)(p-2) - 1} + b(x)\cdot \mathrm{sgn}(x)
    > 0 \qquad\mbox{ and }
\\
& (p-1)^2 \theta^{p-1} (\theta - 1)\,
    |x|^{(\theta - 1)(p-2) - 1} + \theta\, b(x)\cdot \mathrm{sgn}(x)
  - \lambda\, |x|
    \leq 0 \,.
\end{align*}
The last two inequalities, respectively, are equivalent with
\begin{equation}
\label{ineq:b(x)/x}
\begin{aligned}
& {}- (p-1)\left[ p\theta^{p-1} - (p-1)\theta^{p-2}\right]\,
      |x|^{(\theta - 1)(p-2) - 2}
  < b(x) / x
\\
  \leq
& {}- (p-1)^2\, \theta^{p-2}(\theta - 1)\, |x|^{(\theta - 1)(p-2) - 2}
  + \theta^{-1}\, \lambda
    \quad\mbox{ whenever }\, 0 < |x| < 1 \,.
\end{aligned}
\end{equation}
Hence, the validity of the desired inequality \eqref{e:df/d_theta>0}
is derived from our choice of
the {\it\bfseries convection velocity\/} $b(x)$ in \eqref{e:b(x)},
with a help from eq.~\eqref{e:varphi(x)} as follows:

First, the inequality on the left\--hand side of \eqref{ineq:b(x)/x}
is justified by combining the inequalities in \eqref{e:p,theta}
with $0 < |x|\leq 1$ in the calculations that follow:
\begin{align}
\nonumber
&   \frac{b(x)}{x}
  + (p-1)\left[ p\theta^{p-1} - (p-1)\theta^{p-2}\right]\,
    |x|^{(\theta - 1)(p-2) - 2}
\\
\nonumber
&
\begin{aligned}
& {}
  = (p-1)^2 (\theta_1 - 1) \theta_1^{p-2}\,
    |x|^{ (p-2) (\theta_2 - 1) - 2 }
\\
& {}
  + (p-1)\, [ p\theta - (p-1) ]\, \theta^{p-2}\,
    |x|^{(\theta - 1)(p-2) - 2}
\end{aligned}
\\
\label{eq:b(x)}
&
\begin{aligned}
& {}
  \geq
    (p-1)^2 (\theta_1 - 1) \theta_1^{p-2}\,
    |x|^{ (p-2) (\theta_2 - 1) - 2 }
\\
& {}
  + (p-1)\, [ p\theta_1 - (p-1) ]\, \theta_1^{p-2}\,
    |x|^{(\theta_2 - 1)(p-2) - 2}
\end{aligned}
\\
\nonumber
& = (p-1)\, \theta_1^{p-2}\,
    [ (p-1) (\theta_1 - 1) + p\, (\theta_1 - 1) + 1 ]\,
    |x|^{ (p-2) (\theta_2 - 1) - 2 }
\\
\nonumber
& = (p-1)\, \theta_1^{p-2}\, [ (2p-1) (\theta_1 - 1) + 1 ]\,
    |x|^{ (p-2) (\theta_2 - 1) - 2 }
  > 0 \,.
\end{align}

Similarly, the inequality on the right\--hand side of \eqref{ineq:b(x)/x}
is justified by
\begin{equation}
\label{eq:varphi(x)}
\begin{aligned}
&   \frac{b(x)}{x}
  + (p-1)^2\, \theta^{p-2}(\theta - 1)\, |x|^{(\theta - 1)(p-2) - 2}
\\
&
\begin{aligned}
& {}
  = (p-1)^2 (\theta_1 - 1) \theta_1^{p-2}\,
    |x|^{ (p-2) (\theta_2 - 1) - 2 }
\\
& {}
  + (p-1)^2\, \theta^{p-2} (\theta - 1)\, |x|^{(\theta - 1)(p-2) - 2}
\end{aligned}
\\
& < 2 (p-1)^2\, \theta_2^{p-2} (\theta_2 - 1)
  < 2 (p-1)^2\, (\theta_2 - 1) \theta_2^{p-1}\, \theta^{-1}
  \leq \theta^{-1}\, \lambda \,,
\end{aligned}
\end{equation}
owing to
$\lambda\geq 2 (p-1)^2\, (\theta_2 - 1) \theta_2^{p-1} > 0$.

The ``counterexample claims'' in Example~\ref{exam-counter}
thus follow from \eqref{e:df/d_theta>0}.
%
\hfill\Square
\par\vskip 10pt

We {\bf remark\/} that the inequality
``lower bound $\,<\,$ upper bound"
in ineq.~\eqref{ineq:b(x)/x} above is satisfied thanks to
\begin{equation*}
  {}- (p-1)\left[ p\theta^{p-1} - (p-1)\theta^{p-2}\right]\,
  \leq
  {} - (p-1)^2\theta^{p-2}(\theta - 1)
\end{equation*}
combined with
$\lambda = \varphi'\left( 1 - |x|^{\theta}\right) > 0$
whenever $0 < |x| < 1$.
Indeed, we have
\begin{equation}
\label{cond:b(x)/x}
\begin{aligned}
    \left[ p\theta^{p-1} - (p-1)\theta^{p-2}\right]
  - (p-1)\theta^{p-2}(\theta - 1)
& {}
  = (p-1)\theta^{p-2}
    \left( \frac{p}{p-1}\, \theta - 1 - \theta + 1\right)
\\ 
& {}
  = (p-1)\theta^{p-1} \,\frac{1}{p-1} = \theta^{p-1} > 1 \,.
\end{aligned}
\end{equation}
As far as our counterexample is concerned,
this observation allows us to choose any continuous function
$b\colon [-1,1]\to \mathbb{R}$ satisfying ineq.~\eqref{ineq:b(x)/x}
in Example~\ref{exam-counter} above.
Our choice of $b(x)$ in eq.~\eqref{e:b(x)}
is even continuously differentiable and its derivative $b'(x)$
vanishes at zero,
\begin{math}
  b'(0) = \lim_{x\to 0} b(x)/x = 0\,,
\end{math}
owing to the inequalities in \eqref{e:p,theta}.
%
\par\vskip 10pt

Our counterexample (Example~\ref{exam-counter})
allows for the following slight modification / generalization.

\begin{example}[{\rm A modification of Example~\ref{exam-counter}}]
\label{exam-counter_mod}
\begingroup\rm
In addition to eqs.\ \eqref{e:p,theta} and \eqref{e:varphi(x)},
instead of eq.~\eqref{e:b(x)} in Example~\ref{exam-counter} above,
another choice of the convection velocity $b(x)$ can be any
continuously differentiable function
$b\colon [-1,1]\to \mathbb{R}$ satisfying the following conditions:
\begin{align*}
&
\left\{\hspace{2.00mm}
\begin{aligned}
& b(0) = b'(0) = 0 \quad\mbox{ and }\quad
  b'(x)\geq 0 \quad\mbox{ for every }\, x\in [-1,1] \,,
\\
& 0\leq \frac{b(x)}{x}\leq 
        (p-1)^2 (\theta_1 - 1) \theta_1^{p-2}\,
        |x|^{ (p-2) (\theta_2 - 1) - 2 }
      \quad\mbox{ whenever }\, 0 < |x|\leq 1 \,.
\end{aligned} 
\right.
\end{align*}

Then the desired inequalities in \eqref{ineq:b(x)/x}
are obtained in much the same way as in
Example~\ref{exam-counter}, namely, from inequalities very similar to
\eqref{eq:b(x)} and \eqref{eq:varphi(x)}.
\endgroup
%
\hfill\Square
\end{example}
\par\vskip 10pt

We {\bf remark\/} that in both our counterexamples,
Example \ref{exam-counter} and~\ref{exam-counter_mod},
the {\it\bfseries SCP\/} fails to hold precisely at the point $x=0$,
\begin{equation*}
  \{ 0\} = \{ x\in (-1,1)\colon b(x) = 0\} \,,
\end{equation*}
where the convection velocity $b(x)$ vanishes.
In contrast with this observation,
in Example~\ref{exam-SCP-Girg} (in {\S}\ref{ss:strong_crit} below)
we will see another entirely different counterexample to
the {\it\bfseries SCP\/} with
$b(x)\equiv b_0 = \mathrm{const} < 0$ being a negative constant and
$\lambda = 0$.
Moreover, the equality $u(x) = v(x)$ will hold throughout the interval
$[-1,\, -1/2]\subset [-1,1]$ whereas
$u(x) < v(x)$ holds for every $x\in (-1/2,\, 1)$.

It is also necessary to stress that the ``natural sufficient condition''
on the functions $b$ and $\varphi$,
stated in ineq.~\eqref{hypo:b,phi} in Paragraph {\S}\ref{ss:weak-sufficient},
which guarantees the validity of the {\em\bfseries WCP\/},
is trivially fulfilled in both our
Examples \ref{exam-counter} and~\ref{exam-counter_mod}, thanks to
\begin{equation*}
\begin{aligned}
  \frac{1}{2}\, b'(x) +
  \frac{\partial\varphi}{\partial u}(x,s)\geq
  \frac{1}{2}\, b'(x) +
  \inf_{s'\in \RR}
  \frac{\partial\varphi}{\partial u}(x,s') =
  \frac{1}{2}\, b'(x) + \lambda\geq \lambda > 0
\\
  \quad\mbox{ for }\,
    (x,s)\in [-1,1]\times \mathbb{R} \,,
\end{aligned}
\end{equation*}
where
$\varphi\colon [-1,1]\times \mathbb{R}\to \mathbb{R}$
is given by
$\varphi(x,s)\equiv \varphi(s) = \lambda\, s$
for $(x,s)\in [-1,1]\times \mathbb{R}$,
by eq.~\eqref{e:varphi(x)}.
Under this condition, ineq.~\eqref{hypo:b,phi},
the {\em\bfseries WCP\/} will be established
in Proposition~\ref{prop-Weak_Comp} below.
%
\par\vskip 10pt

\subsection{Interpretation for non\--Newtonian fluid flows
            {\rm (e.g., of a pollutant)}}
\label{ss:intro-fluid}

The case of vanishing {\it convection velocity\/}
$\mathbf{b}\colon \Omega\to \mathbb{R}^N$
in the general $N$\--dimensional elliptic problem~\eqref{e:ell_PDE}
leads to a very popular variational problem for
the {\it\bfseries critical points\/}
of the corresponding energy functional
\begin{equation}
\label{energy:ell_PDE}
\begin{aligned}
& \mathcal{E}(u)\eqdef
    \frac{1}{p}\int_{\Omega} |\nabla u|^p \,\mathrm{d}x
  + \int_{\Omega} \Phi(x,u) \,\mathrm{d}x \,,\quad
    u\in W_0^{1,p}(\Omega) \,;
\\
&   \quad\mbox{ with }\quad
  \Phi(x,s)\eqdef \int_0^s \varphi(x,t) \,\mathrm{d}t
    \quad\mbox{ for }\,
    (x,s)\in \Omega\times \mathbb{R} \,,
\end{aligned}
\end{equation}
defined on the Sobolev space $W_0^{1,p}(\Omega)$.
This variational problem has been investigated in numerous articles
during at least the last seven decades.
A classical reference for this quasilinear elliptic problem is
the monograph by {\sc J.-L.\ Lions} \cite{J-L_Lions1969}.
More recent results can be found in
{\sc P.\ Dr\'abek} \cite{P_Drabek1992},
{\sc P.\ Dr\'abek}, {\sc A.\ Kufner}, and {\sc F.\ Nicolosi}
\cite{Drabek-Kuf-Nic-1997}, and
{\sc P.\ Tak\'a\v{c}} \cite{P_Takac2004},
and in the numerous references therein
(\cite{P_Drabek1992, Drabek-Kuf-Nic-1997, P_Takac2004}).

The variational approach does not allow for
the {\it\bfseries convection term\/}
$\mathbf{b}(x)\cdot \nabla u$ in problem~\eqref{e:ell_PDE}.
It can be treated by methods of 
{\it\bfseries pseudo\--monotone operators\/}
(\cite[Chapt.~2, {\S}2.4 -- {\S}2.5, pp.\ 179--182]{J-L_Lions1969})
or fixed points for the variational problem
(\cite[Chapt.~2, {\S}2.6, pp.\ 182--190]{J-L_Lions1969}).
A more recent approach to the non\--variational problem~\eqref{e:ell_PDE}
with a convection term can be found in
{\sc J.\ Garc\'{\i}a\--Meli\'an}, {\sc J.~C.\ Sabina de Lis}, and
{\sc P.\ Tak\'a\v{c}} \cite{Garcia-Lis-Takac2017}.

Finally, if the WCP is valid
(as stated in Proposition~\ref{prop-Weak_Comp} below),
then also {\it\bfseries monotone methods\/}
can be applied to the non\--variational problem~\eqref{e:ell_PDE}.

\begin{quote}
From the {\rm engineer's\/} point of view
(cf.\ \cite{Benedikt-Kotrla2018}),
the {\it\bfseries WCP\/}
means that increasing the content of certain substance
(e.g., a pollutant such as oil or nitrate)
in the model described by the elliptic partial differential equation
\eqref{e:ell_PDE} by means of the source function $f(x)$
leads to an increase of the substance density $u(x)$
throughout the entire domain  $\Omega\subset\mathbb{R}^N$.
The point here is that the increase from $u(x)$ to, say,
$v(x)$ does not have to be strict at every point $x\in \Omega$,
that is, the strict inequality $u(x) < v(x)$ might not hold in spite of
the increase of $f(x)$ to $g(x)$ with $f(x)\leq g(x)$
for every $x\in \Omega$ and $f\not\equiv g$ in $\Omega$.
Here, we have replaced the pair $(f,u)$ by the new pair $(g,v)$
in eq.~\eqref{e:ell_PDE}; cf.\ eqs.~\eqref{eq:WCP}.

In particular, if pollution by the substance with
the density $u(x)$ in the model is caused locally by
the source function $f(x)$,
by increasing the source function $f(x)$ to $g(x)$ with $f(x)\leq g(x)$
and $f\not\equiv g$,
a global pullution increase
from the substance density $u(x)$ to $v(x)$ with $u(x) < v(x)$
(throughout the entire domain $\Omega$),
occurs if and only if the {\it\bfseries SCP\/} is valid.
\end{quote}

As usual, given a pair of functions
$f,g\in L^{\infty}(\Omega)$,
the relation $f\not\equiv g$ in $\Omega$ means that the set
\begin{math}
  \{ x\in \Omega\colon f(x)\neq g(x)\}
\end{math}
has positive $N$\--dimensional Lebesgue measure in $\RR^N$.
If $f,g\in C(\Omega)$, then
$f\not\equiv g$ in $\Omega$ means that the set
\begin{math}
  \{ x\in \Omega\colon f(x)\neq g(x)\}
\end{math}
contains an open ball\hfil\break
\begin{math}
  B_r(x_0)\eqdef \{ x\in \mathbb{R}^N\colon |x-x_0| < r\}
          \subset \Omega
\end{math}
of radius $r > 0$.

\section{The Weak Comparison Principle}
\label{s:weak_comp}

In this section we focus on the {\em\bfseries WCP\/},
a weaker form of the {\em\bfseries Comparison Principle\/}
that we wish to investigate.
To this end, let us consider
the two\--point boundary value problems \eqref{eq:WCP}
stated in the Introduction
(Section~\ref{s:intro}, {\S}\ref{ss:weak-counter}).
In the next paragraph, {\S}\ref{ss:weak-sufficient},
we provide a simple, natural sufficient condition on the functions
$b\colon [-1,1]\to \mathbb{R}$ and
$\varphi\colon [-1,1]\times \mathbb{R}\to \mathbb{R}$
which guarantee the validity of the {\em\bfseries WCP\/},
that is to say,
$u(x)\leq v(x)$ for all $x\in [-1,1]$.

\subsection{A Natural Sufficient Condition}
\label{ss:weak-sufficient}

We assume that the functions $b(x)$, $\varphi(x,u)$, $f(x)$, and $g(x)$
satisfy the following hypothesis:

\par\vskip 10pt
%
\begin{itemize}
\item[{\bf (M)}]
$b\colon [-1,1]\to \mathbb{R}$ is continuously differentiable,
both
\begin{math}
  \varphi ,\; \frac{\partial\varphi}{\partial u}
  \colon [-1,1]\times \mathbb{R}\to \mathbb{R}
\end{math}
are continuous, and the following
{\em\bfseries ``monotonicity''\/} inequality is satisfied:
\begin{equation}
\label{hypo:b,phi}
  \frac{1}{2}\, b'(x) +
  \frac{\partial\varphi}{\partial u}(x,s)\geq 0
  \quad\mbox{ for all }\, (x,s)\in (-1,1)\times \mathbb{R} \,.
\end{equation}
In addition, we assume that $f,g\in L^{\infty}(-1,1)$.
\end{itemize} 
%
\par\vskip 20pt

In order to obtain a weak comparison principle for problem
\eqref{eq:WCP} analogous to Theorem 10.7 in
\cite[Sect.~10, {\S}10.4, p.~268]{GilbargTrudinger2001},
we establish the following rather general result:

\par\vskip 10pt
\begin{proposition}[{\rm Weak Comparison Principle}]
\label{prop-Weak_Comp}
Assume\/ {\rm Hypothesis~{\bf (M)}}.
If\/ $f,g \in L^{\infty}(-1,1)$ satisfy\/
$f\leq g$ a.e.\ in $(-1,1)$, and\/
$u(-1)\leq v(-1)$ and\/ $u(+1)\leq v(+1)$, then
$u\leq v$ holds everywhere in $[-1,1]$.
\end{proposition}
\par\vskip 10pt

{\it Proof.\/}
We follow similar ideas as does the proof of Theorem 10.7
(under condition {\rm (iii)})
in \cite[pp.\ 270--271]{GilbargTrudinger2001}.
Subtracting the differential equations for $u$ and $v$ in~\eqref{eq:WCP}
from each other, we arrive at
\begin{equation*}
\begin{aligned}
& {}
  - \frac{\mathrm{d}}{\mathrm{d}x}
    \left( |u'|^{p-2} u' - |v'|^{p-2} v'\right)
  - b(x)\, (u'-v')
  + \left[ \varphi(x,u) - \varphi(x,v)\right]
\\
& {}
  = f(x) - g(x)\leq 0 \quad\mbox{ for }\, x\in (-1,1) \,.
\end{aligned}
\end{equation*}
We multiply this inequality by the function
$w\stackrel{\mathrm{def}}{=} (u-v)^+$,
$w\in W_0^{1,p}(-1,1)$,
which has the vanishing boundary values $w(\pm 1) = 0$,
and integrate with respect to $x\in (-1,1)$, thus arriving at
\begin{equation}
\begin{aligned}
\label{e:(u-v)^+}
& \int_{-1}^1
    \left( |u'|^{p-2} u' - |v'|^{p-2} v'\right) w'(x) \,\mathrm{d}x
  - \int_{-1}^1 b(x)\, (u'-v')\, w(x) \,\mathrm{d}x
\\
& {}
  + \int_{-1}^1
    \left[ \varphi(x,u) - \varphi(x,v)\right] w(x) \,\mathrm{d}x
  = \int_{-1}^1 (f-g)\, w(x) \,\mathrm{d}x\leq 0 \,.
\end{aligned}
\end{equation}
Since the boundary conditions
$u(-1)\leq v(-1)$ and $u(+1)\leq v(+1)$ guarantee
$w(\pm 1) = 0$, the second integral on the left\--hand side of
eq.~\eqref{e:(u-v)^+} becomes
\begin{equation}
\begin{aligned}
\label{e:b(u'-v')(u-v)^+}
& \int_{-1}^1 b(x)\, (u'-v')\, w(x) \,\mathrm{d}x
  = \frac{1}{2}\int_{-1}^1
    b(x)\, \frac{\mathrm{d}}{\mathrm{d}x}\, [w(x)]^2 \,\mathrm{d}x
  = {}- \frac{1}{2}\int_{-1}^1
        b'(x)\, [w(x)]^2 \,\mathrm{d}x \,.
\end{aligned}
\end{equation}
In order to treat the first and third integrals on the left\--hand side
of eq.~\eqref{e:(u-v)^+}, we use the standard formulas
\begin{align}
\label{e:(u'-v')}
    |u'|^{p-2} u' - |v'|^{p-2} v'
& {}
  = (p-1)\left[ \int_0^1
  \left| (1-\theta)\, u'(x) + \theta\, v'(x)\right|^{p-2} \,\mathrm{d}\theta
         \right] (u'-v') \,,
\\
\label{e:phi(u)-phi(v)}
    \varphi(x,u) - \varphi(x,v)
& {}
  = \left[ \int_0^1 \frac{\partial\varphi}{\partial u}
    \bigl( x,\, (1-\theta)\, u(x) + \theta\, v(x)\bigr) \,\mathrm{d}\theta
    \right] (u-v) \,.
\end{align}
We insert eqs.\ \eqref{e:b(u'-v')(u-v)^+}, \eqref{e:(u'-v')}, and
\eqref{e:phi(u)-phi(v)} into eq.~\eqref{e:(u-v)^+}, thus obtaining
\begin{equation*}
\begin{aligned}
& (p-1)\int_{-1}^1\left[ \int_0^1
  \left| (1-\theta)\, u'(x) + \theta\, v'(x)\right|^{p-2} \,\mathrm{d}\theta
         \right] [w'(x)]^2 \,\mathrm{d}x
  + \frac{1}{2}\int_{-1}^1
        b'(x)\, [w(x)]^2 \,\mathrm{d}x
\\
& {}
  + \int_{-1}^1\left[ \int_0^1 \frac{\partial\varphi}{\partial u}
    \bigl( x,\, (1-\theta)\, u(x) + \theta\, v(x)\bigr) \,\mathrm{d}\theta
               \right] [w(x)]^2 \,\mathrm{d}x
  = \int_{-1}^1 (f-g)\, w(x) \,\mathrm{d}x\leq 0
\end{aligned}
\end{equation*}
or, equivalently,
\begin{equation}
\begin{aligned}
\label{ineq:(u-v)^+}
& (p-1)\int_{-1}^1\left[ \int_0^1
  \left| (1-\theta)\, u'(x) + \theta\, v'(x)\right|^{p-2} \,\mathrm{d}\theta
         \right] [w'(x)]^2 \,\mathrm{d}x
\\
& {}
  + \int_{-1}^1\left[ \int_0^1
    \left( \frac{1}{2}\, b'(x) + \frac{\partial\varphi}{\partial u}
    \Bigl( x,\, (1-\theta)\, u(x) + \theta\, v(x)\Bigr)
    \right) \,\mathrm{d}\theta
               \right] [w(x)]^2 \,\mathrm{d}x
\\
& {}
  = \int_{-1}^1 (f-g)\, w(x) \,\mathrm{d}x\leq 0 \,.
\end{aligned}
\end{equation}

Finally, we take $s = (1-\theta)\, u(x) + \theta\, v(x)$
in inequality \eqref{hypo:b,phi} and insert it into \eqref{ineq:(u-v)^+};
it forces $w'(x) = 0$ almost everywhere throughout $(-1,1)$
as the first integral on the left\--hand side must be non\--positive.
The boundary conditions $w(\pm 1) = 0$ then yield
$w(x) = (u(x) - v(x))^{+} = 0$ for all $x\in [-1,1]$, that is,
we have proved that $u\leq v$ in $[-1,1]$, as claimed.
\hfill\rule{3mm}{3mm}
\par\vskip 10pt

Although the following two corollaries are concerned with
an {\it\bfseries ``interior''\--type\/}
{\it\bfseries SCP\/}, they are direct consequences of
the {\it\bfseries WCP\/}
that we have just proved above in Proposition~\ref{prop-Weak_Comp}.
This is the reason why we have decided to state them right after
our proof of Proposition~\ref{prop-Weak_Comp}.
In both corollaries that follow next, we assume that the reaction function
$\varphi(x,s)\equiv \varphi(x)$ is independent from the dummy variable
$s\in \mathbb{R}$.
In this case, the monotonicity inequality \eqref{hypo:b,phi}
is reduced to $b'(x)\geq 0$ for all $x\in (-1,1)$.
We may take $\varphi(x,s)\equiv 0$ without loss of generality by replacing
the pair of functions $f(x)$ and $g(x)$ in eqs.~\eqref{eq:WCP} by
$f(x) - \varphi(x)$ and $g(x) - \varphi(x)$,
respectively, thus arriving at
\begin{equation}
\label{eq:WCP:varphi=0}
\left\{\quad
\begin{aligned}
&
{}- \left( |u'|^{p-2} u'\right)' - b(x)\, u' = f(x) \,,
\qquad
{}- \left( |v'|^{p-2} v'\right)' - b(x)\, v' = g(x) \,;
\\
& \quad\mbox{ with }\quad
  u(\pm 1)\leq v(\pm 1) \quad\mbox{ and }\quad
  f(x)\leq g(x) \;\mbox{ for }\, x\in (-1,1) \,.
\end{aligned}
\right.
\end{equation}
%

\begin{corollary}[{\rm ``Interior'' Strong Comparison Principle}]
\label{cor-Int_Strong_Comp}
In the situation of\/ {\rm Proposition~\ref{prop-Weak_Comp}\/},
with $\varphi(x, \,\cdot\,)\equiv 0$,
we have the following two (mutually exclusive) alternatives:
\begin{itemize}
\item[{\bf (a)}]
The set\/
\begin{align*}
  P_1\eqdef \{ x\in (-1,1)\colon u(x) < v(x)\}
\end{align*}
is {\em not empty\/}.
Then there is a pair of points $-1 < a_{-1} < a_1 < 1$ such that\/
$u(a_{-1}) < v(a_{-1})$ and\/ $u(a_1) < v(a_1)$.
Given any such a pair, set\/
\begin{equation*}
  \eta = \min\{ v(a_{-1}) - u(a_{-1}) ,\, v(a_1) - u(a_1) \} > 0 \,.
\end{equation*}
Furthermore, in problem~\eqref{eq:WCP:varphi=0} above we have\/
$v(x)\geq u(x) + \eta > u(x)$ for every $x\in [a_{-1},\, a_1]$.
Finally, $P_1$ is an open interval in $\mathbb{R}$.
\item[{\bf (b)}]
The set\/ $P_1$ is {\em empty\/}.
This is the case if and only if all\/
$u(-1) = v(-1)$, $u(1) = v(1)$, and $f = g$ a.e.\ in $(-1,1)$
hold in the boundary value problems in \eqref{eq:WCP:varphi=0}.
\end{itemize}
\end{corollary}
\par\vskip 10pt

{\it Proof of\/} Corollary~\ref{cor-Int_Strong_Comp}.
Let us begin with Alternative {\bf (a)}, i.e., let
$P_1\not= \emptyset$.
Thanks to the continuity of both functions
$u,v\colon [-1,1]\to \mathbb{R}$,
$P_1$ is a nonempty open subset of the interval $(-1,1)$.
Consequently, there are two points
$-1 < a_{-1} < a_1 < 1$ such that
\begin{equation*}
  \eta = \min\{ v(a_{-1}) - u(a_{-1}) ,\, v(a_1) - u(a_1) \} > 0 \,.
\end{equation*}
We apply the {\it\bfseries WCP\/}
from Proposition~\ref{prop-Weak_Comp} on the interval
$[a_{-1},\, a_1]$ in place of $[-1,1]$ and for the pair of functions
$u(x) + \eta$, $v$ instead of $u(x)$, $v(x)$, respectively.
By these arguments, the inequalities
$u(a_{-1}) + \eta\leq v(a_{-1})$ and
$u(a_1) + \eta\leq v(a_1)$ imply $u(x) + \eta\leq v(x)$
for every $x\in [a_{-1},\, a_1]$.
Consequently, $P_1$ is an open and connected subset of $(-1,1)$
and, thus, an open interval itself.

In contrast, Alternative {\bf (b)} forces $u\equiv v$
throughout the entire interval $[-1,1]$.
Problems \eqref{eq:WCP:varphi=0} then yield $f = g$ a.e.\ in $(-1,1)$.

Our corollary is proved.
\hfill\rule{3mm}{3mm}
\par\vskip 10pt

\section{The Strong Comparison Principle near a Critical Point}
\label{s:strong_critical}

In this section we focus on
some rather technical issues and open questions related to
the {\it\bfseries SCP\/}
for problems \eqref{eq:WCP:varphi=0}
(that is to say, for problems \eqref{eq:WCP}
 with the vanishing reaction function $\varphi(x,s)\equiv 0$)
in case when the two solutions,
$u,v\colon [-1,1]\to \mathbb{R}$, respectively,
satisfy the weak comparison principle and they meet in an interior point
$x_0\in (-1,1)$, that is,
$u(x)\leq v(x)$ for all $x\in [-1,1]$ and $u(x_0) = v(x_0)$.
With regard to
Corollary~\ref{cor-Int_Strong_Comp}, Alt.~{\bf (b)}, above,
we will focus our attention on the most problematic case when
\begin{itemize}
\item[{\rm (i)}]
$u(x) = v(x)$ holds for every $x\in [-1,x_0]$ whereas
$u(x) < v(x)$ for every $x\in (x_0,a_1)$
\hfil\break
in a nonempty open interval $(x_0,a_1)$, for some $a_1\in (x_0,1]$,
i.e., $P_1 = (x_0,a_1)\neq \emptyset$.
\end{itemize}
The symmetric ``mirror'' case below can be treated analogously,
\begin{itemize}
\item[{\rm (i')}]
$u(x) = v(x)$ holds for every $x\in [x_0,1]$ whereas
$u(x) < v(x)$ for every $x\in (a_{-1},x_0)$
\hfil\break
in a nonempty open interval $(a_{-1},x_0)$, for some $a_{-1}\in [-1,x_0)$,
i.e., $P_1 = (a_{-1},x_0)\neq \emptyset$.
\end{itemize}

Throughout this entire section
{\em we focus on the interaction\/} between
the {\em\bfseries diffusion\/} and {\em\bfseries convection\/} effects,
that is, between the terms
$\left( |u'|^{p-2} u'\right)'$ and $b(x)\, u'$
$\bigl(
 \left( |v'|^{p-2} v'\right)'$ and $b(x)\, v'$,
respectively$\bigr)$
in eqs.~\eqref{eq:WCP:varphi=0},
while completely ignoring the dependence on the dummy variable
$s\in \mathbb{R}$ in the reaction function $\varphi(x,s)$
that stands for the value $s = u(x)$ ($s = v(x)$, respectively).
In other words, we take
$\varphi(x,s)\equiv \varphi(x)$ independent from $s\in \mathbb{R}$
and include it in the ``source'' function $f(x)$
on the right\--hand side of eq.~\eqref{eq:intro:counter}
by replacing $f(x)$ by the difference $f(x) - \varphi(x)$.
This reduction hypothesis thus allows us to simply set
$\varphi(x,s)\equiv 0$ in eq.~\eqref{eq:intro:counter}; cf.\
problem~\eqref{eq:WCP:varphi=0} above.

\subsection{The solutions $u$ and $v$ near a Common Interior Meeting Point}
\label{ss:strong_meeting}

We begin by recalling a special case from the work by
{\sc M.\ Cuesta} and {\sc P.\ Tak\'{a}\v{c}}
\cite[Sect.~4, Example 4.1, pp.\ 740--741]{CuestaTakac2000}
that we have adapted to our model in eq.~\eqref{eq:intro:counter}
in Example~\ref{exam-counter} above.
In this counterexample to the SCP we have $x_0 = 0$.

Recalling our notation from Corollary~\ref{cor-Int_Strong_Comp},
we {\bf denote} by $P_1$ and $P_0$ the set of all points
$x\in (-1,1)$ where the {\it\bfseries SCP\/}
is valid, i.e., $u(x) < v(x)$, and is {\bf not\/} valid, i.e.,
$u(x) = v(x)$, respectively:
\begin{align*}
  P_1\eqdef \{ x\in (-1,1)\colon u(x) < v(x)\} \quad\mbox{ and }\quad
  P_0\eqdef \{ x\in (-1,1)\colon u(x) = v(x)\} \,;
\end{align*}
hence,
\begin{math}
  P_1\cup P_0 = \{ x\in (-1,1)\colon u(x)\leq v(x)\} = (-1,1) ,
\end{math}
by our hypothesis on the validity of the {\it\bfseries WCP\/}.
Clearly, $P_1\cap P_0 = \emptyset$.
From this {\it\bfseries WCP\/}, that is,
$u(x)\leq v(x)$ for every $x\in [-1,1]$, we derive easily
\begin{equation*}
  x_0\in P_0\subset (-1,1) \;\Longrightarrow\; u'(x_0) = v'(x_0)
  \quad\mbox{ in addition to }\quad u(x_0) = v(x_0) \,.
\end{equation*}
As usual, we denote by $\overline{P}_1$ the closure in
$[-1,1]\subset \mathbb{R}$ of the set $P_1$.

In the remaining part of this paragraph, {\S}\ref{ss:strong_meeting},
we will investigate the behavior of the two functions, $u$ and~$v$,
in the vicinity of a point $x_0\in P_0\cap \overline{P}_1$
with the property \eqref{e:x_0+-delta} below:

\par\vskip 10pt
\begin{lemma}[{\rm Weak Comparison Principle at} $x_0$]
\label{lem-Weak_Comp:x_0}
Assume that\/ $x_0\in P_0\cap \overline{P}_1$ satisfies the following\/

\begin{quote}
{\bf Hypothesis\/}:
There is a number $\delta > 0$ such that
\begin{equation}
\label{e:x_0+-delta}
  (x_0 - \delta, x_0)\subset P_1 \quad\mbox{ or }\quad
  (x_0, x_0 + \delta)\subset P_1 \,.
\end{equation}
%
\end{quote}

\par\noindent
Then we have\/ $u(x_0) = v(x_0)$ and\/ $u'(x_0) = v'(x_0) = 0$.
\end{lemma}
\par\vskip 10pt

We postpone the proof until after Remark~\ref{rem-Weak_Comp:x_0}.

\par\vskip 10pt
\begin{corollary}[{\rm {\it\bfseries WCP\/} at} $x_0$]
\label{cor-WCP:f(x_0)<=g(x_0)}
Let\/ $x_0\in P_0\cap \overline{P}_1$ satisfy\/
{\rm Hypothesis~(\ref{e:x_0+-delta})} and assume that both functions\/
$f,g\colon (-1,1)\subset \mathbb{R}\to \mathbb{R}$
are continuous at the point $x_0$.
Then we have also\/
\begin{equation}
\label{e:Weak_Comp:x_0}
  \frac{\mathrm{d}}{\mathrm{d}x}\,
  \left( |u'|^{p-2} u'\right) \Big\vert_{x = x_0}
  = {}- f(x_0)\geq {}- g(x_0) =
  \frac{\mathrm{d}}{\mathrm{d}x}\,
  \left( |v'|^{p-2} v'\right) \Big\vert_{x = x_0} \,.
\end{equation}
\end{corollary}
\par\vskip 10pt

Recalling that
$b\colon [-1,1]\to \mathbb{R}$ is continuously differentiable,
by our {\rm Hypothesis~{\bf (M)}} with Ineq.~\eqref{hypo:b,phi},
from Corollary~\ref{cor-WCP:f(x_0)<=g(x_0)} above
we are able to derive the following stronger result:

\par\vskip 10pt
\begin{corollary}[{\rm Vanishing Convection/Drift at} $x_0$]
\label{cor-WCP:f(x_0)=g(x_0)}
In the situation of\/ {\rm Corollary~\ref{cor-WCP:f(x_0)<=g(x_0)}}
above, with $x_0\in P_0\cap \overline{P}_1$ satisfying\/
{\rm Hypothesis~(\ref{e:x_0+-delta})} and with both functions\/
$f,g\colon (-1,1)\subset \mathbb{R}\to \mathbb{R}$
being continuous at the point $x_0$, we must have the equality\/
\begin{equation}
\label{eq:Weak_Comp:x_0}
  {}- \frac{\mathrm{d}}{\mathrm{d}x}\,
  \left( |u'|^{p-2} u'\right) \Big\vert_{x = x_0}
  = f(x_0) = g(x_0) =
  {}- \frac{\mathrm{d}}{\mathrm{d}x}\,
  \left( |v'|^{p-2} v'\right) \Big\vert_{x = x_0} \,.
\end{equation}
\end{corollary}
\par\vskip 10pt

The following remark on {\bf Hypothesis~(\ref{e:x_0+-delta})}
is {\em\bfseries important\/} for the results
in Lemma~\ref{lem-Weak_Comp:x_0}
and Corollaries \ref{cor-WCP:f(x_0)<=g(x_0)}
and~\ref{cor-WCP:f(x_0)=g(x_0)}.
Namely, if the following stronger condition
\begin{equation*}
  (x_0 - \delta, x_0)\subset P_1 \quad\mbox{ and }\quad
  (x_0, x_0 + \delta)\subset P_1
\end{equation*}
is valid, then we must have $u(x_0) < v(x_0)$ as well, by
Corollary~\ref{cor-Int_Strong_Comp}, Alt.~{\bf (a)}
(cf.\ also Corollary~\ref{cor-Int_Strong_Comp}, Alt.~{\bf (b)}).

\begin{remark}\label{rem-Weak_Comp:x_0}\nopagebreak
\begingroup\rm
Assume that\/
\begin{equation}
\label{eq:x_0+-delta}
  (x_0 - \delta, x_0)\cup (x_0, x_0 + \delta)\subset P_1 \,.
\end{equation}
Then also $u(x_0) < v(x_0)$ holds, by
Alt.~{\bf (a)} of Corollary~\ref{cor-Int_Strong_Comp}, i.e.,
$x_0\not\in P_0$ and, consequently, we get
\begin{math}
  (x_0 - \delta, x_0 + \delta)\subset P_1 \,.
\end{math}
\hfill\Square
\endgroup
\end{remark}
\par\vskip 10pt

\begin{quote}
In contrast, our counterexample to the {\it\bfseries SCP\/}
(Example~\ref{exam-counter})
in the Introduction (Section~\ref{s:intro}), {\S}\ref{ss:intro-math},
shows that the equality $u(x_0) = v(x_0)$ is still possible
in the more general model governed by
the boundary value problem \eqref{eq:intro:counter}
which allows also for the reaction function $\varphi(u)$
taken to be equal to
$\varphi(x,u) = \lambda\, |u|^{p-2} u$
in \cite[Example 4.1]{CuestaTakac2000} and to
$\varphi(x,u) = \lambda\, u$
in Examples \ref{exam-counter} and~\ref{exam-counter_mod} above.

From the {\rm engineer's\/} point of view
(cf.\ \cite{Benedikt-Kotrla2018}),
in this counterexample
the {\bf reaction\/} of the reaction function $\varphi(u)$
is so strong near the point $x_0 = 0$ that at the point $x_0 = 0$
it {\bf eliminates\/} an increase of the source function $f(x)$
throughout the two neighboring intervals in condition
\eqref{eq:x_0+-delta} above.
Thus, the equality $u(x_0) = v(x_0)$ holds in this counterexample
at $x_0 = 0$.

In particular, if pollution by the substance with
the density $u(x)$ in the model is caused globally in the region
\begin{equation*}
  (-1,0)\cup (0,1) = (-1,1)\setminus \{ 0\} \,,
\end{equation*}
increasing the source function $f(x)$ to, say, $g(x)$ globally
throughout this entire region, $(-1,1)\setminus \{ 0\}$,
except for the point $x_0 = 0$, i.e.,
$f(x) < g(x)$ whenever $0 < |x| < 1$ and $f(0) = g(0)$,
{\bf does not mean\/} that an increase of pollution must be felt also
at the point $x_0 = 0$.
\end{quote}

\par\vskip 10pt
{\it\bfseries Proof of\/} Lemma~\ref{lem-Weak_Comp:x_0}.
The two unknown functions
$u,v\colon (-1,1)\subset \mathbb{R}\to \mathbb{R}$
satisfy the following two\--point boundary value problems on $(-1,1)$,
respectively; cf.\ \eqref{eq:WCP:varphi=0}:
\begin{equation}
\label{eqs:SCP}
\left\{\hspace{2.00mm}
\begin{aligned}
&
{}- \left( |u'|^{p-2} u'\right)' - b(x) u' = f(x) \,,
\qquad
{}- \left( |v'|^{p-2} v'\right)' - b(x) v' = g(x) \,;
\\
& u(\pm 1) = v(\pm 1) = 0 \,, \quad\mbox{ and } \quad
  f\leq g \;\mbox{ a.e.\ in }\, (-1,1) \,,
  \quad f\not\equiv g \;\mbox{ in } (-1,1) \,.
\end{aligned}
\right.
\end{equation}

By contradiction to the claim $u'(x_0) = v'(x_0) = 0$,
suppose that either
$u'(x_0) = v'(x_0) > 0$ or else $u'(x_0) = v'(x_0) < 0$.
Subtracting the latter from the former equation in problems
\eqref{eqs:SCP} above, we arrive at
the linear ``elliptic'' differential equation
\begin{equation}
\label{e:thm:lineariz}
\left\{\hspace{5.00mm}
\begin{aligned}
  {}- [D(u',v')(x) w'(x)]' - b(x)\, w'(x) & = g(x) - f(x)\ ({}\geq 0)
    \quad\mbox{ in }\, (-1,1) \,;
\\
  w(\pm 1) & = 0 \,,
\end{aligned}
\right.
\end{equation}
for the difference
$w = v-u\colon (-1,1)\subset \mathbb{R}\to \mathbb{R}$ on $(-1,1)$.
This equation is {\em regular\/} at and near the point
$x_0\in (-1,1)$, meaning that the diffusion coefficient
$D(u',v')(x)$ given by
\begin{equation}
\label{e:D(u',v')}
\begin{aligned}
  D(u',v')(x)
& {}
  = (p-1)\int_0^1
    \left| (1-\theta) u'(x) + \theta v'(x)\right|^{p-2} \,\mathrm{d}\theta
\\
& {}
  \geq c\left(
        \max_{0\leq \theta\leq 1} |(1-\theta) u'(x) + \theta v'(x)|
        \right)^{p-2}
  \geq 0 \quad\mbox{ on }\, [-1,1]
\end{aligned}
\end{equation}
satisfies $D(u',v')(x) > 0$ for every
$x\in (-1,1)$ near the point $x_0\in (-1,1)$.
Here, $c > 0$ is some constant independent of $u'(x)$, $v'(x)$,
cf.\ \cite[Ineq.\ (A.1), p.~233]{Takac2002}.
In general, problem \eqref{e:thm:lineariz}
is not uniformly elliptic.
However, since the inequality $D(u',v')(x) > 0$ holds strictly
provided at least one of the two inequalities
$u'(x)\neq 0$ and $v'(x)\neq 0$ is valid,
the problem for the function $w = v - u$ is, in fact, strictly elliptic.

Since also $x_0\in \overline{P}_1\cap (-1,1)$ and there is a number
$\delta > 0$ such that {\rm Hypothesis~(\ref{e:x_0+-delta})} is valid,
we may apply the boundary point principle
({\em Hopf's lemma\/})
established recently in the article by
{\sc D.~E.\ Apushkinskaya} and {\sc A.~I.\ Nazarov}
\cite[Theorem 2.1, p.~681]{Apush-Nazarov-2019}.
Owing to $w(x) = v(x) - u(x) > 0$ for all
$x\in (x_0 - \delta, x_0)$ or for all
$x\in (x_0, x_0 + \delta)$, and to $w(x_0) = v(x_0) - u(x_0) = 0$,
we thus conclude that $w'(x_0) < 0$.
This inequality clearly contradicts the obvious inequality
(the {\it\bfseries WCP\/} established in
 Proposition~\ref{prop-Weak_Comp})
$w(x) = v(x) - u(x)\geq 0$ for all points $x\in (-1,1)$ near $x_0$,
where $w(x_0) = v(x_0) - u(x_0) = 0$.

We have proved that $u'(x_0) = v'(x_0) = 0$.
\hfill\rule{3mm}{3mm}
\par\vskip 10pt

{\it\bfseries Proof of\/} Corollary~\ref{cor-WCP:f(x_0)<=g(x_0)}.
We insert the conclusion of Lemma~\ref{lem-Weak_Comp:x_0},
particularly $u'(x_0) = v'(x_0) = 0$, into eqs.~\eqref{eqs:SCP}
in order to obtain ineq.~\eqref{e:Weak_Comp:x_0}
in Corollary~\ref{cor-WCP:f(x_0)<=g(x_0)}.
\hfill\rule{3mm}{3mm}
\par\vskip 10pt

{\it\bfseries Proof of\/} Corollary~\ref{cor-WCP:f(x_0)=g(x_0)}.
Let us suppose that, by contradiction
to eq.~\eqref{eq:Weak_Comp:x_0}, we have the sharp inequality
in ineq.~\eqref{e:Weak_Comp:x_0}, that is,
\begin{math}
  {}- f(x_0) > {}- g(x_0) .
\end{math}
Set
\begin{equation*}
  \eta\eqdef \genfrac{}{}{}1{1}{2} \left( {}- f(x_0) + g(x_0)\right) > 0 \,.
\end{equation*}
All three functions
$b, f, g\colon (-1,1)\subset \mathbb{R}\to \mathbb{R}$
being continuous at the point $x_0$, we conclude from
ineq.~\eqref{e:Weak_Comp:x_0} that there is a number
$\delta_1\in (0,\delta)$ such that the inequality
\begin{equation}
\label{diff:Weak_Comp:x_0}
\begin{aligned}
  \frac{\mathrm{d}}{\mathrm{d}x}\,
  \left( |u'|^{p-2} u' - |v'|^{p-2} v'\right) \Big\vert_{x}
  = {}- b(x)\, [u'(x) - v'(x)] - f(x) + g(x)
  \geq \eta > 0
\\
    \quad\mbox{ holds for every }\,
    x\in (x_0 - \delta_1,\, x_0 + \delta_1) \,.
\end{aligned}
\end{equation}
We first recall that $u'(x_0) = v'(x_0) = 0$,
by Lemma~\ref{lem-Weak_Comp:x_0}, then integrate the inequality in
\eqref{diff:Weak_Comp:x_0} above from $x_0$ to
$x\in (x_0 - \delta_1,\, x_0 + \delta_1)$, thus arriving at
\begin{alignat}{2}
\label{int:Weak_Comp:x<x_0}
  \left( |u'|^{p-2} u' - |v'|^{p-2} v'\right) \Big\vert_{x}
& \leq {}- \eta\, |x-x_0|
&&  \quad\mbox{ for every }\, x\in (x_0 - \delta_1,\, x_0] \,,
\\
\label{int:Weak_Comp:x>x_0}
  \left( |u'|^{p-2} u' - |v'|^{p-2} v'\right) \Big\vert_{x}
& \geq \eta\, |x-x_0|
&&  \quad\mbox{ for every }\, x\in [x_0,\, x_0 + \delta_1) \,.
\end{alignat}
Consequently, we have $u'(x) - v'(x) < 0$
throughout the interval $(x_0 - \delta_1,\, x_0)$.
When combined with $u(x_0) = v(x_0)$ and $u'(x_0) = v'(x_0)$,
this inequality forces
\begin{equation*}
\begin{aligned}
  u(x) - v(x) = [u(x) - v(x)] - [u(x_0) - v(x_0)]
  = [u'(\xi) - v'(\xi)] (x - x_0) > 0
\\
    \quad\mbox{ for every }\, x\in (x_0 - \delta_1,\, x_0)
    \,\mbox{ and some }\, \xi\in (x,x_0)\subset (x_0 - \delta_1,\, x_0) \,,
\end{aligned}
\end{equation*}
by the mean value theorem.
But already this inequality, i.e., $u(x) - v(x) > 0$ for all
$x\in (x_0 - \delta_1,\, x_0)$,
contradicts the weak comparison principle $u\leq v$ in $[-1,1]$.

Alternatively to ineq.~\eqref{int:Weak_Comp:x<x_0},
we may take advantage of the inequality in
\eqref{int:Weak_Comp:x>x_0} in a similar way: it yields
$u'(x) - v'(x) > 0$ on the interval $(x_0,\, x_0 + \delta_1)$.
When combined with $u(x_0) = v(x_0)$ and $u'(x_0) = v'(x_0)$,
this inequality forces
\begin{equation*}
\begin{aligned}
  u(x) - v(x) = [u(x) - v(x)] - [u(x_0) - v(x_0)]
  = [u'(\xi) - v'(\xi)] (x - x_0) > 0
\\
    \quad\mbox{ for every }\, x\in (x_0,\, x_0 + \delta_1)
    \,\mbox{ and some }\, \xi\in (x_0,x)\subset (x_0,\, x_0 + \delta_1) \,,
\end{aligned}
\end{equation*}
by the mean value theorem.
Thus, we get $u(x) - v(x) > 0$ for all $x\in (x_0,\, x_0 + \delta_1)$
which contradicts the weak comparison principle
$u\leq v$ in $[-1,1]$ again.

The corollary is proved.
\null\hfill\rule{3mm}{3mm}
\par\vskip 10pt

The following proposition is a slight extension of
Corollary~\ref{cor-Int_Strong_Comp}
({\rm ``Interior'' Strong Comparison Principle}):

\par\vskip 10pt
\begin{proposition}[{\rm ``Interior'' Strong Comparison Principle}]
\label{prop-Interval_SCP}
In the situation of\/ {\rm Proposition~\ref{prop-Weak_Comp}\/}
and\/ {\rm Corollary~\ref{cor-Int_Strong_Comp}\/},
with $\varphi(x, \,\cdot\,)\equiv 0$, the set\/
\begin{align*}
  P_1 = \{ x\in (-1,1)\colon u(x) < v(x)\}
  \quad\mbox{ is an open interval in $\mathbb{R}$. }
\end{align*}
If it is empty, then $u\equiv v$ and\/ $f\equiv g$ hold throughout
the interval\/
$P_0 = (-1,1)$.
If\/ $P_1\neq \emptyset$, then we have\/
$P_1 = (a_{-1},\, a_1)$ with a pair of points\/
$-1\leq a_{-1} < a_1\leq 1$;
moreover, only the following three possibilities may occur:
\begin{itemize}
\item[{\rm (i)}]$\;$
$P_1 = (-1,1)$ in which case the
{\em strong comparison principle\/} ({\it\bfseries SCP\/})
is valid:
$u(x) < v(x)$ holds for all\/ $x\in (-1,1)$.
\item[{\rm (ii)}]$\;$
$-1 < a_{-1} < a_1\leqq 1$ in which case we have\/
$u(x) = v(x)$ and $f(x) = g(x)$ for all\/ $x\in (-1,\, a_{-1}]$,
together with the vanishing derivatives\/
$u'(a_{-1}) = v'(a_{-1}) = 0$.
\item[{\rm (iii)}]$\;$
$-1\leqq a_{-1} < a_1 < 1$ in which case we have\/
$u(x) = v(x)$ and $f(x) = g(x)$ for all\/ $x\in [a_1,\, 1)$,
together with the vanishing derivatives\/
$u'(a_1) = v'(a_1) = 0$.
\end{itemize}
\end{proposition}
\par\vskip 10pt

{\it Proof of\/} Proposition~\ref{prop-Interval_SCP}.
Corollary~\ref{cor-Int_Strong_Comp} guarantees that the set $P_1$
is an open interval in $\mathbb{R}$, possibly empty.
If $P_1 = \emptyset$, then obviously $P_0 = (-1,1)$.
The complementary alternative, $P_1\neq \emptyset$, forces
$P_1 = (a_{-1},\, a_1)$ with a pair of points
$-1\leqq a_{-1} < a_1\leqq 1$, by
Alt.~{\bf (b)} of Corollary~\ref{cor-Int_Strong_Comp}.
Case~{\rm (i)} is obvious.

If Case~{\rm (i)} is false, i.e., $P_0\neq \emptyset$,
then at least one of the two
Cases {\rm (ii)} and {\rm (iii)} must occur:$\;$
The conclusion in both these cases follows from
Corollary~\ref{cor-Int_Strong_Comp}, Alt.~{\bf (a)}, combined with
Lemma~\ref{lem-Weak_Comp:x_0}.
\null\hfill\rule{3mm}{3mm}
\par\vskip 10pt

\subsection{The solutions $u$ and $v$ near a Common Critical Meeting Point}
\label{ss:strong_crit}

In this paragraph we ``refine'' the counterexample to
the {\it\bfseries SCP\/} presented in Example~\ref{exam-counter}
in the Introduction (Section~\ref{s:intro}), {\S}\ref{ss:weak-counter}.
We will construct a pair of functions
$f,g\in L^{\infty}(-1,1)$ such that
\begin{equation*}
  f\leq g \;\mbox{ a.e.\ in }\, (-1,1) \quad\mbox{ and }\quad
  f\not\equiv g \;\mbox{ in }\, (-1,1) \,,
\end{equation*}
and the corresponding pair of solutions
$u,v\colon [-1,1]\to \mathbb{R}$
to the two\--point boundary value problems \eqref{eqs:SCP}
stated in the previous paragraph, {\S}\ref{ss:strong_meeting},
such that in {\rm Hypothesis~(\ref{e:x_0+-delta})} we have
\begin{equation}
\label{e:x_0<x_0+delta}
  P_0 = (-1, x_0] \quad\mbox{ and }\quad
  P_1 = (x_0, 1) \quad\mbox{ with }\,
  x_0 = {}- \genfrac{}{}{}1{1}{2} \in (-1,1) \,.
\end{equation}

\par\vskip 10pt
\begin{example}[{\rm Another Counterexample to the} {\it\bfseries SCP\/}]
\label{exam-SCP-Girg}
\begingroup\rm
In this counterexample we continue using the notation introduced in
{\rm Proposition~\ref{prop-Interval_SCP}}
from the previous paragraph ({\S}\ref{ss:strong_crit}).
Here, we will construct an example of {\rm Part~(ii)}
in Proposition~\ref{prop-Interval_SCP} with the pair of points
$a_{-1} = {}- 1/2 < a_1 = 1$, that is to say,
\begin{itemize}
\item[{\bf (ii)}]$\;$
$-1 < a_{-1} = {}- 1/2 < a_1 = 1$ in which case we have\/
$u(x) = v(x)$ and $f(x) = g(x)$ for all\/ $x\in (-1,\, a_{-1}]$,
together with the vanishing derivatives\/
$u'(a_{-1}) = v'(a_{-1}) = 0$, and
$u(x) < v(x)$ for all\/ $x\in (a_{-1}, 1)$ together with
$u(\pm 1) = v(\pm 1) = 0$.
\end{itemize}

We feel that our construction will be more understandable
when we start from the graphs of the solutions
$u,v\colon [-1,1]\to \mathbb{R}$
to the two\--point boundary value problems \eqref{eqs:SCP}
sketched below in {\bf Figure~1}
where we have introduced the abbreviation
\begin{equation*}
  x_{\mathrm{s}}\eqdef 1 - 3 / 2^{\frac{p}{p - 1}} = 1 - 3 / 2^{p'}
  \in \left( - \genfrac{}{}{}1{1}{2} \,,\, 0\right)
  \quad\mbox{ with the conjugate exponent }\;
  p'= \genfrac{}{}{}1{p}{p-1} \in (1,2) \,.
\end{equation*}
%

%
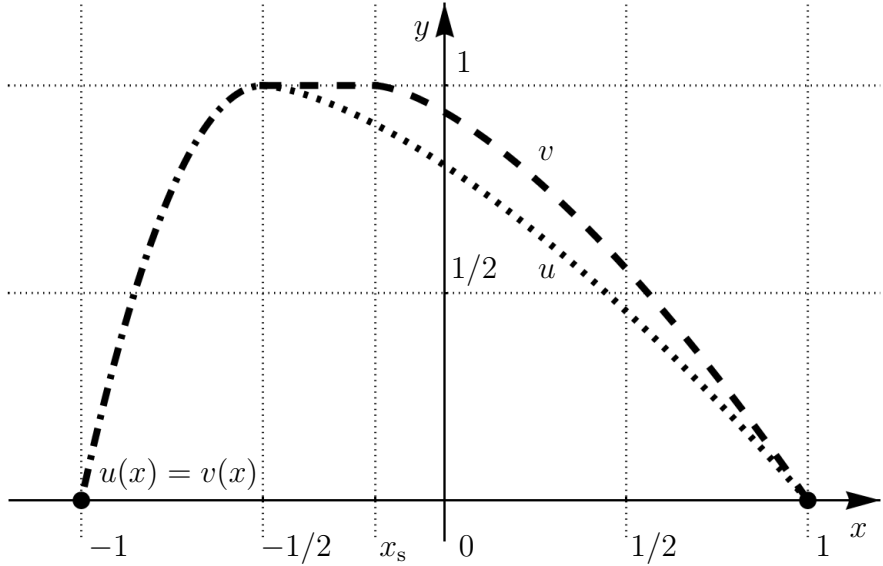
\begin{figure}[ht]
\begin{picture}(200,200)(0,0)
\put(-70,0){\includegraphics[height=7.2cm]
{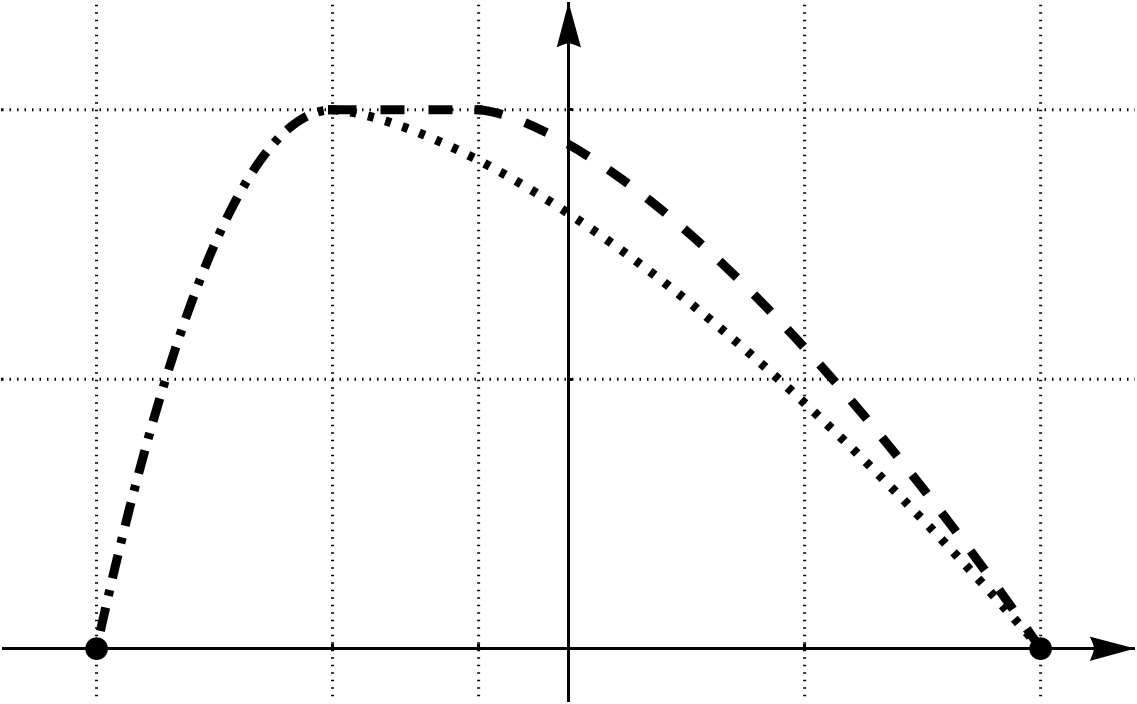}
\put(-130, 100){$u$}
\put(-130,145){$v$}
\put(-296,23){$u(x)=v(x)$}
\put(-300,-5){$-1$}
\put(-235,-5){$-1/2$}
\put(-163,100){$1/2$}
\put(-161,178){$1$}
\put(-190,-5){$x_\mathrm{s}$}
\put(-160,-5){$0$}
\put(-95,-5){$1/2$}
\put(-25,-5){$1$}
\put(-12,2){$x$}
\put(-177,192){$y$}
}
\end{picture}
\label{fig:counterexample}
\caption{Counterexample to the strong comparison principle for $p=4$. 
Here, $a_{-1} = - 1/2$, $x_{\mathrm{s}} = 1 - 3/2^{\frac{4}{3}}$
and the solutions $u$ and $v$ are given by the respective formulas
\eqref{def_u} and \eqref{def_v}.
}
\end{figure}
%
\par\vskip 10pt

Indeed, we will show that the functions
\begin{align}
\label{def_u}
& u(x) \stackrel{\rm def}{=} \left\{
        \begin{aligned}
                & {}- 4x(x+1)\,,
                &&\quad x \in \left[ -1 \,,\, -\frac{1}{2}\right)\,,
\\[5pt]
                & 1 - 3^{\frac{1 - p}{p - 2}}(2x + 1)^{\frac{p - 1}{p - 2}}\,, 	
                &&\quad x \in \left[ -\frac{1}{2} \,,\, 1\right]\,,
        \end{aligned}
\right.
\end{align}
and
\begin{align}
\label{def_v}
& v(x) \stackrel{\rm def}{=} \left\{
        \begin{aligned}
                & {}- 4x(x+1)\,,
                &&\quad x \in \left[ -1 \,,\, -\frac{1}{2}\right)\,,
\\[5pt]
                &\qquad 1 \,,
                &&\quad x \in
  \left[ -\frac{1}{2} \,,\, x_{\mathrm{s}} \right] \,,
\\[5pt]
 		&1 - 2^{\frac{p}{p - 2}}\cdot 3^{\frac{1 - p}{p - 2}}
  \cdot \left( x - x_{\mathrm{s}} \right)^{\frac{p-1}{p-2}} \,,
                &&\quad x \in
        \left( x_{\mathrm{s}} \,,\, 1 \right] \,,
        \end{aligned}
\right.
\end{align}
satisfy \eqref{eqs:SCP} with $p > 2$,
\begin{align*}
  b = - \genfrac{(}{)}{}0{3}{2}^{1 - p}
        \genfrac{(}{)}{}0{p - 2}{p - 1}^{1 - p} < 0
    \quad\mbox{ (constant {\it convection velocity\/}), }
\end{align*}
\begin{align*}
& f(x) \stackrel{\rm def}{=}
\left\{
        \begin{aligned}
 		&2^{p-1}
  \left( 2^p (p-1) [ {}- (2x+1) ]^{p-2}
       - 12\genfrac{(}{)}{}0{p-2}{p-1}^{1-p}\, (2x+1)
  \right)\,, \hskip 5pt
 		&& x\in \left( -1 \,,\, -\frac{1}{2}\right)\,,
\\[5pt]
 		&\qquad 0 \,, 	
 		&& x\in \left[ -\frac{1}{2} \,,\, 1\right)\,,
        \end{aligned}
\right.
\end{align*}
and
\begin{align*}
& g(x) \stackrel{\rm def}{=}
\left\{
        \begin{aligned}
 		&2^{p-1}
  \left( 2^p (p-1) [ {}- (2x+1) ]^{p-2}
       - 12\genfrac{(}{)}{}0{p-2}{p-1}^{1-p}\, (2x+1)
  \right)\,, \hskip 5pt
 		&& x\in \left(-1 \,,\, - \frac{1}{2}\right) \,,
\\[5pt]
 		&\qquad 0 \,,
 		&& x\in
  \left[ -{} \frac{1}{2}\,,\, x_{\mathrm{s}} \right]\,,
\\[5pt]
                &2^{\frac{(p-1)^2 + 1}{p-2}}
                \cdot 3^{{}- \frac{(p-1)^2}{p-2}}
                \cdot \genfrac{(}{)}{}0{p-1}{p-2}^p 
  \left( x - x_{\mathrm{s}} \right)^{\frac{1}{p-2}} \,,\,
                && x\in
                \left( x_{\mathrm{s}} \,,\, 1\right) \,.
        \end{aligned}
\right.
\end{align*}
Note that, thanks to $p > 2$, both functions $f$ and $g$
are bounded on $[-1,1]$, and thus $f,g\in L^{\infty}(-1,1)$.
It is easily seen that
$f(x) = g(x)\geq 0$ in $\left( -1,\, - 1/2 \right)$,
$f(x) = g(x) = 0$ in
$\left[ -1/2 ,\, x_{\mathrm{s}} \right]$, and
$g(x) > f(x) = 0$ in $(x_{\mathrm{s}}, 1)$.
We conlude that $f(x)\leq g(x)$ in $(-1,1)$.

By our definition of functions $u$, $v$, $f$, and $g$,
we easily verify the following equalities,
\begin{equation*}
  f(x) = {}- \left( |u'(x)|^{p-2} u'(x)\right)' - b\, u'(x)
    \quad\mbox{ and }\quad
  g(x) = {}- \left( |v'(x)|^{p-2} v'(x)\right)' - b\, v'(x) \,.
\end{equation*}
The reader can verify these two relations by investigating
the particular cases of $x\in (-1,1)$ suggested by
the partition of the interval $(-1,1)$ sketched
in {\bf Figure~1} above, that is to say,
%
\begin{equation}
\label{e:partition}
  (-1,1) = \left( -1,\, - 1/2 \right)
    \cup   \left[ - 1/2 ,\, x_{\mathrm{s}} \right]
    \cup   \left( x_{\mathrm{s}} ,\, 0 \right]
    \cup   (0,1) \,.
\end{equation}
This is a matter of a direct calculation in each of the subintervals
indicated above.

By our definition of functions $u$ and $v$, we have
$u\equiv v$ in $(-1,\, -1/2)$.
The strict inequality $u(x) < v(x)$ for every $x\in (-1/2 ,\, 1)$
is verified again by investigating the various subintervals of
the partition \eqref{e:partition}; cf.\ {\bf Figure~1}.
\endgroup
%
\hfill\Square
\end{example}
\par\vskip 10pt

In order to prevent the equality $u(x)\equiv v(x)$ on the interval
$(-1,\, -1/2)$ in our Example~\ref{exam-SCP-Girg} above,
we make the following natural hypothesis:

\par\vskip 20pt
\par\noindent
{\bf Hypothesis~$\mathbf{H}_{\pm 1}$.}
Let $f,g\in L^{\infty}(-1,1)$ and
$0\leq f(x)\leq g(x)$ hold for almost all $x\in (-1,1)$.
In addition, we assume that for every $\delta > 0$ small enough we have
$f\not\equiv g$ in both intervals
$(-1, -1+\delta)$ and $(1-\delta, 1)$, that is to say, both sets
\begin{equation*}
  \{ x\in (-1,1)\colon f(x) < g(x)\} \cap (-1, -1+\delta)
  \quad\mbox{ and }\quad
  \{ x\in (-1,1)\colon f(x) < g(x)\} \cap (1-\delta, 1)
\end{equation*}
have positive Lebesgue measure for every $\delta > 0$ small enough.
\par\vskip 20pt

Consequently, {\rm Hypothesis}~$\mathbf{H}_{\pm 1}$
is valid with any number
$\delta\in \left( 0 ,\, 1/2\right)$
in Lemma~\ref{lem-Hopf_Max:+-1} below.
This lemma will then yield
$v'(1) < u'(1)\leq 0\leq u'(-1) < v'(-1)$.

As an immediate consequence we obtain
the following boundary point lemma of Hopf type.

\begin{lemma}[{\rm Hopf's Maximum Principle at} $x = \pm 1$]
\label{lem-Hopf_Max:+-1}
Assume that\/ {\rm Hypothesis}~$\mathbf{H}_{\pm 1}$ is valid.
Then the following {\em\bfseries boundary point lemma\/} holds:\/
$v'(-1) < 0 < v'(1)$.

In particular, the linear elliptic equation
\eqref{e:thm:lineariz} for the difference
$w\stackrel{\mathrm{def}}{=} v - u$ on $(-1,1)$
is regular in $(-1,1)$ near the endpoints $x = \pm 1$, i.e.,
there is some $\delta > 0$ small enough such that
the inequality $D(u',v')(x)\ge \mathrm{const} > 0$
holds strictly for every\/
\begin{math}
  x\in (-1, -1+\delta)\cup (1-\delta, 1) \,.
\end{math}

Consequently, we have even
the \emph{Hopf boundary point comparison\/}:
\hfil\break
$v'(1) < u'(1)\leq 0\leq u'(-1) < v'(-1)$.
\end{lemma}
\par\vskip 10pt

{\it Proof of\/} Lemma~\ref{lem-Hopf_Max:+-1}.
We apply the {\em\bfseries Hopf boundary point lemma\/} from
{\sc P.\ Tolksdorf}
\cite[Prop.\ 3.2.1 and 3.2.2, p.~801]{Tolksdorf-1983} or
{\sc J.~L.\ V\'azquez} \cite[Theorem~5, p.~200]{Vazquez-1984}.
This result, due to {\sc E.~Hopf}, is described and proved in all details
in the monograph by
{\sc P.\ Pucci} and {\sc J.\ Serrin}
\cite[{\S}5.4]{Pucci-Serr_book}, Theorem 5.5.1 on p~120; cf.\ also
\cite[{\S}5.1]{Pucci-Serr_book}, Example on p.~104.

Owing to $g\geq 0$ in $(-1,1)$ and
$g\not\equiv 0$ in both intervals
$(-1, -1+\delta)$ and $(1-\delta, 1)$,
for every $\delta > 0$ small enough, we thus get $v'(-1) < 0 < v'(1)$.
In particular, the linear elliptic equation
\eqref{e:thm:lineariz} for the difference
$w\stackrel{\mathrm{def}}{=} v - u$ on $(-1,1)$ being regular,
we arrive at
$w'(-1) < 0 < w'(1)$ by the new version of
the {\em Hopf boundary point lemma\/}
established recently in the article by
{\sc D.~E.\ Apushkinskaya} and {\sc A.~I.\ Nazarov}
\cite{Apush-Nazarov-2019}.

Hence, the desired inequalities
$v'(1) < u'(1)\leq 0\leq u'(-1) < v'(-1)$
follow.
\hfill\rule{3mm}{3mm}
\par\vskip 10pt

\begin{theorem}[{\rm Hopf's ``Global'' Strong Comparison Principle}]
\label{thm-Hopf's_SCP}
Let us consider\/ {\rm problem\/}~\eqref{eqs:SCP}
for $u$ and $v$ on $(-1,1)$.
If also\/
$f,g\in L^{\infty}(-1,1)$ with\/
$0\leq f\leq g$ a.e.\ in $(-1,1)$ and\/
{\rm Hypothesis}~$\mathbf{H}_{\pm 1}$ is valid, then
{\em\bfseries Hopf's (``Global'') Strong Comparison Principle}
holds:
\begin{equation*}
  v'(1) < u'(1)\leq 0\leq u'(-1) < v'(-1)
  \quad\mbox{ and }\quad
  0\leq u(x) < v(x) \quad\mbox{ for all }\, x\in (-1,1) \,.
\end{equation*}
\end{theorem}
\par\vskip 10pt

{\it Proof of\/} Theorem~\ref{thm-Hopf's_SCP}.
Let us recall that both functions
$u,v\colon [-1,1]\to \mathbb{R}$ are continuously differentiable.
We begin by applying Proposition~\ref{prop-Interval_SCP}
to conclude that the set
\begin{equation*}
  P_1 = \{ x\in (-1,1)\colon u(x) < v(x)\}
  \quad\mbox{ is an open interval in $\mathbb{R}$, i.e., }\
\end{equation*}
$P_1 = (a_{-1},\, a_1)$ holds with a pair of points
$-1\leq a_{-1} < a_1\leq 1$; hence, the complementary set in $(-1,1)$
equals to
\begin{align*}
  P_0 = \{ x\in (-1,1)\colon u(x) = v(x)\}
  = (-1,a_{-1}]\cup [a_1,1) \,.
\end{align*}
On the other hand, by Lemma~\ref{lem-Hopf_Max:+-1} above,
if we choose $\delta > 0$ small enough, such that both inequalities,
$u'(\xi) < v'(\xi)$ and $v'(\eta) < u'(\eta)$, are valid for all
$\xi\in (-1, -1+\delta)$ and $\eta\in (1-\delta, 1)$,
then we must have also
$u(\xi) < v(\xi)$ and $u(\eta) < v(\eta)$ again for all
$\xi\in (-1, -1+\delta)$ and $\eta\in (1-\delta, 1)$,
thanks to the Dirichlet boundary conditions
$u(-1) = v(-1)$ and $u(1) = v(1)$.
But this result forces
$(-1, -1+\delta)\cup (1-\delta, 1)\subset P_1$
which clearly contradicts the statement
\begin{equation*}
  (-1,1)\setminus P_1 = P_0 = (-1,a_{-1}]\cup [a_1,1)
  \quad\mbox{ unless }\, P_0 = \emptyset \,,
\end{equation*}
i.e., $P_1 = (-1,1)$ holds as desired.
Lemma~\ref{lem-Hopf_Max:+-1} yields also
\begin{math}
  v'(1) < u'(1)\leq 0\leq u'(-1) < v'(-1) \,.
\end{math}
The remaining part of {\rm Hopf's Comparison Principle}
follows directly from $P_1 = (-1,1)$.

The theorem is proved.
\hfill\rule{3mm}{3mm}
\par\vskip 10pt

\subsection{The case of constant convection velocity}
\label{ss:constant_b}

In this paragraph we will assume that the convection velocity
$b\colon (-1,1)\to \mathbb{R}$ is a given constant, i.e.,
$b(x)\equiv b_0\in \mathbb{R}$ for all $x\in (-1,1)$.
Hence, Proposition~\ref{prop-Weak_Comp}
applies and, thus, the {\it\bfseries WCP\/} is valid.
In this simplified model, we will be able to prove a weaker form of
Theorem~\ref{thm-Hopf's_SCP},
Theorem~\ref{thm-b_0-Hopf's_SCP} below,
under the following hypothesis,

\par\vskip 10pt
\par\noindent
{\bf Hypothesis~$\mathbf{H}_0$.}
Let $f,g\in L^{\infty}(-1,1)$ and
$0\leq f(x)\leq g(x)$ hold for almost all $x\in (-1,1)$
together with $f\not\equiv g$ in the interval $(-1,1)$, i.e., the set
\begin{equation*}
  \{ x\in (-1,1)\colon f(x) < g(x)\}
\end{equation*}
has positive Lebesgue measure.
In addition, we assume that
$\varphi(x, \,\cdot\,)\equiv 0$ and
the {\em\bfseries convection velocity\/}
$b(x) = b_0\equiv \mathrm{const}$ is constant.
\par\vskip 10pt

We note that $\mathbf{H}_0$ is weaker than
{\rm Hypothesis}~$\mathbf{H}_{\pm 1}$
used in Lemma~\ref{lem-Hopf_Max:+-1} and Theorem~\ref{thm-Hopf's_SCP}.
This weaker result, in Theorem~\ref{thm-b_0-Hopf's_SCP} below,
is motivated and closely related to our results above,
in Paragraph {\S}\ref{ss:strong_meeting}.

\par\vskip 10pt
\begin{theorem}[{\rm Another Hopf's ``Global'' Comparison Principle}]
\label{thm-b_0-Hopf's_SCP}
Let us consider\/ {\rm problem\/}~\eqref{eqs:SCP}
for $u$ and $v$ on $(-1,1)$, with\/
$b(x) = b_0\equiv \mathrm{const}$ being constant on $(-1,1)$.
Assume that\/
$f,g\in L^{\infty}(-1,1)$ satisfy\/
{\rm Hypothesis}~$\mathbf{H}_0$ and\/
$u(x_0) = v(x_0)$ holds at some point\/
$x_0\in (-1,1)$.
Then we have also\/
$u'(x_0) = v'(x_0) = 0$ together with the following two alternatives:
\begin{itemize}
\item[{\bf Alt.~1:}]
If\/ $b_0\leq 0$ then we must have
$u(x) = v(x)$ for all $ x\in (-1,x_0]$ together with\/
$f(x) = g(x)$ for almost every $x\in (-1,x_0)$, that is, both,
$u\equiv v$ and $f\equiv g$ in $(-1,x_0)$.
\item[{\bf Alt.~2:}]
If\/ $b_0\geq 0$ then we must have
$u(x) = v(x)$ for all $ x\in [x_0,1)$ together with\/
$f(x) = g(x)$ for almost every $x\in (x_0,1)$, that is, both,
$u\equiv v$ and $f\equiv g$ in $(x_0,1)$.
\end{itemize}
\end{theorem}
\par\vskip 10pt

\begin{remark}\label{rem-b_0-Hopf's_SCP}\nopagebreak
\begingroup\rm
The result in {\bf Alt.~1} reflects the scenario in
{\rm Example~\ref{exam-SCP-Girg}} with\/
\begin{math}
  b_0 = - \genfrac{(}{)}{}1{3}{2}^{1 - p}
          \genfrac{(}{)}{}1{p - 2}{p - 1}^{1 - p}\hfil\break < 0
\end{math}
and\/ $x_0 = -1/2$.
\hfill\Square
\endgroup
\end{remark}
\par\vskip 10pt

{\it\bfseries Proof of\/} Theorem~\ref{thm-b_0-Hopf's_SCP}.
Let us recall that both functions
$u,v\colon [-1,1]\to \mathbb{R}$ are continuously differentiable.
We begin by applying Proposition~\ref{prop-Interval_SCP}
to conclude that the set
\begin{equation*}
  P_1 = \{ x\in (-1,1)\colon u(x) < v(x)\}
  \quad\mbox{ is an open interval in $\mathbb{R}$, i.e., }\
              P_1 = (a_{-1},\, a_1)
\end{equation*}
with a pair of points $-1\leq a_{-1} < a_1\leq 1$.
Hence, its complementary set in $(-1,1)$ equals to
\begin{align*}
  P_0 = \{ x\in (-1,1)\colon u(x) = v(x)\}
  = (-1,a_{-1}]\cup [a_1,1) \,.
\end{align*}

We rewrite the general {\rm problem\/}~\eqref{eqs:SCP}
in the special form with $b(x) = b_0\equiv \mathrm{const}$
being a real constant, $b_0\in \mathbb{R}$,
\begin{equation}
\label{eqs:SCP:b_0}
\left\{\quad
\begin{aligned}
&
{}- \frac{\mathrm{d}}{\mathrm{d}x}\,
    \left( |u'|^{p-2} u' + b_0\, u\right) = f(x) \,,
\qquad
{}- \frac{\mathrm{d}}{\mathrm{d}x}\,
    \left( |v'|^{p-2} v' + b_0\, v\right) = g(x) \,;
\\
& \quad
  u(\pm 1) = v(\pm 1) = 0 \,,
  \quad\mbox{ and }\quad
\\
& \quad
  0\leq f\leq g \quad\mbox{ a.e.\ in }\, (-1,1) \,,
  \quad f\not\equiv g \quad\mbox{ in } (-1,1) \,.
\end{aligned}
\right.
\end{equation}
In analogy with our proof of Theorem~\ref{thm-Hopf's_SCP} above,
we wish to show that either
$P_1 = (x_0,\, 1)$ or else $P_1 = (-1,\, x_0)$ holds,
this time without assuming the validity of
{\rm Hypothesis}~$\mathbf{H}_{\pm 1}$.
The desired equivalent result states that either
$P_0 = (-1,\, x_0]$ or else $P_0 = [x_0,\, 1)$, respectively.

Thus, let $x_0\in (-1,1)$ be a point satisfying $u(x_0) = v(x_0)$.
Then the ``local'' inequality $u(x)\leq v(x)$ for all
$x\in (x_0 - \delta,\, x_0 + \delta )\subset (-1,1)$
forces also $u'(x_0) = v'(x_0)$.
Since $u\not\equiv v$ in $(-1,1)$, we can choose $x_0\in (-1,1)$
above in such a way that
{\em at least one\/} of the following two statements is valid:
\begin{quote}
{\rm Ineq.~1:}$\;$
$u(x) < v(x)$ holds for every $x\in (x_0 - \delta,\, x_0)$.
\hfil\break
{\rm Ineq.~2:}$\;$
$u(x) < v(x)$ holds for every $x\in (x_0,\, x_0 + \delta )$.
\end{quote}
(We note that these inequalities may hold simultaneously.)

Furthermore, the linear elliptic equation \eqref{e:thm:lineariz}
for the difference
$w\stackrel{\mathrm{def}}{=} v-u$ on $(-1,1)$
is uniformly elliptic unless $u'(x_0) = v'(x_0) = 0$ is valid.
However, if $u'(x_0) = v'(x_0)\neq 0$,
then the Hopf boundary point lemma would force either
$u'(x_0) > v'(x_0)$ or else $u'(x_0) < v'(x_0)$, respectively.
This would contradict the equality $u'(x_0) = v'(x_0)$.
Thus, we have a point $x_0\in (-1,1)$ such that
$u(x_0) = v(x_0)$ and $u'(x_0) = v'(x_0) = 0$; cf.\
Lemma~\ref{lem-Weak_Comp:x_0}
({\rm {\it\bfseries WCP\/} at} $x_0$).

Next, we use $x_0\in (-1,1)$ as the initial point
for our {\sl ``shooting method''\/}
applied to the special {\rm problem\/}~\eqref{eqs:SCP:b_0} above,
with the initial values $u(x_0) = v(x_0)$ and $u'(x_0) = v'(x_0) = 0$:
\begin{equation*}
\left\{\quad
\begin{aligned}
  |u'(x)|^{p-2} u'(x) + b_0\, u(x) & =
  |u'(x_0)|^{p-2} u'(x_0) + b_0\, u(x_0)
\textstyle
  - \int_{x_0}^x f(s) \,\mathrm{d}s
\\
& = b_0\, u(x_0)
\textstyle
  - \int_{x_0}^x f(s) \,\mathrm{d}s
  \quad\mbox{ for all }\, x\in (-1,1) \,,
\\
  |v'(x)|^{p-2} v'(x) + b_0\, v(x) & =
  |v'(x_0)|^{p-2} v'(x_0) + b_0\, v(x_0)
\textstyle
  - \int_{x_0}^x g(s) \,\mathrm{d}s
\\
& = b_0\, v(x_0)
\textstyle
  - \int_{x_0}^x g(s) \,\mathrm{d}s
  \quad\mbox{ for all }\, x\in (-1,1) \,;
\\
& u(\pm 1) = v(\pm 1) = 0 \,,
  \quad\mbox{ and }\quad
\\
  0\leq f\leq g
& \quad\mbox{ a.e.\ in }\, (-1,1) \,,
  \quad f\not\equiv g \quad\mbox{ in } (-1,1) \,.
\end{aligned}
\right.
\end{equation*}
Subtracting the first differential equation (for $u$)
from the second one (for $v$) we arrive at
\begin{equation}
\label{diff:SCP:b_0}
\left\{\quad
\begin{aligned}
& |v'(x)|^{p-2} v'(x) - |u'(x)|^{p-2} u'(x)
  = {}- b_0 (v(x) - u(x))
\\
&
\textstyle
  {}- \int_{x_0}^x (g(s) - f(s)) \,\mathrm{d}s
  \quad\mbox{ for all }\, x\in (-1,1) \,,
\\
& u(\pm 1) = v(\pm 1) = 0 \,,
  \quad\mbox{ and }\quad
\\
& 0\leq f\leq g
  \quad\mbox{ a.e.\ in }\, (-1,1) \,,
  \quad f\not\equiv g \quad\mbox{ in } (-1,1) \,.
\end{aligned}
\right.
\end{equation}

\underline{\bf Case}~$b_0\leq 0$:
We have
${}- b_0 (v(x) - u(x))\geq 0$ and thus
eq.~\eqref{diff:SCP:b_0} yields
\begin{equation*}
\begin{aligned}
& |v'(x)|^{p-2} v'(x) - |u'(x)|^{p-2} u'(x)
\textstyle
  \geq {}- \int_{x_0}^x (g(s) - f(s)) \,\mathrm{d}s
\\
&
\textstyle
  = \int_x^{x_0} (g(s) - f(s)) \,\mathrm{d}s\geq 0
  \quad\mbox{ for all }\, x\in (-1,x_0] \,.
\end{aligned}
\end{equation*}
Consequently, we must have
$(v-u)'(x)\geq 0$ for all $x\in (-1,x_0]$ which forces
$v(x) - u(x) = 0$ for all $x\in (-1,x_0]$, thanks to
$v(-1) - u(-1) = 0$ and $v(x_0) - u(x_0) = 0$.
This is possible only if $f(x) = g(x)$ holds for almost every
$x\in (-1,x_0)$,
as described in alternative ``{\bf Alt.~1}''.

\underline{\bf Case}~$b_0\geq 0$:
We have
${}- b_0 (v(x) - u(x))\leq 0$ and thus
eq.~\eqref{diff:SCP:b_0} yields
\begin{equation*}
  |v'(x)|^{p-2} v'(x) - |u'(x)|^{p-2} u'(x)
\textstyle
  \leq {}- \int_{x_0}^x (g(s) - f(s)) \,\mathrm{d}s\leq 0
  \quad\mbox{ for all }\, x\in [x_0,1) \,.
\end{equation*}
Consequently, we have
$(v-u)'(x)\leq 0$ for all $x\in [x_0,1)$ whence
$v(x) - u(x) = 0$ for all $x\in [x_0,1)$, thanks to
$v(x_0) - u(x_0) = 0$ and $v(1) - u(1) = 0$.
This entails $f(x) = g(x)$ for almost every $x\in (x_0,1)$,
as claimed in alternative ``{\bf Alt.~2}''.

The proof is finished.
\hfill\rule{3mm}{3mm}
\par\vskip 10pt

{\bf Acknowledgments:}$\;$
The research of
{\sc Ji\v{r}\'{\i} Benedikt}, {\sc Petr Girg}, and {\sc Luk\'a\v{s} Kotrla}
was partially supported by
the {\sl\bfseries Grant Agency of the Czech Republic (GA{\v{C}}R)\/}
under Project No.~{\sf\bfseries 18--03253S\/}.
Ing.\ {\sc Luk\'a\v{s} Kotrla} was partially supported also by
the Ministry of Education, Youth, and Sports
({\sl\bfseries M{\v{S}}MT\/}, Czech Republic)
under the program NPU~I, Project No.~{\sf\bfseries LO1506 (PU--NTIS)\/}.

\par\vskip 20pt

\end{document}